\documentclass[11pt,a4paper]{article}
\usepackage[latin1]{inputenc}
\usepackage{amsmath}
\usepackage{amsfonts}
\usepackage{amssymb}

\newtheorem{satz}{Satz}[section]
\newtheorem{theorem}[satz]{Theorem}
\newtheorem{lemma}[satz]{Lemma}
\newtheorem{cor}[satz]{Corollary}
\newtheorem{defin}[satz]{Definition}
\newtheorem{examples}[satz]{Examples}

\newcommand{\abs}[1]{\left|{#1}\right|}
\newcommand{\rund}[1]{\left(#1\right)}
\newcommand{\spitz}[1]{\left\langle{#1}\right\rangle}
\newcommand{\eckig}[1]{\left[{#1}\right]}
\newcommand{\schweif}[1]{\left\{#1\right\}}
\newcommand{\ceil}[1]{\left\lceil{#1}\right\rceil}

\def\zz{\mathbb{Z}}
\def\cz{\mathbb{C}}
\def\nz{{\rm I\kern-.20em N}}
\def\rz{{\rm I\kern-.20em R}}
\def\pz{{\rm I\kern-.20em P}}

\def\C{\mathcal{C}}
\def\E{\mathcal{E}}
\def\M{\mathcal{M}}
\def\N{\mathcal{N}}
\def\O{\mathcal{O}}
\def\F{\mathcal{F}}
\def\X{\mathcal{X}}
\def\S{\mathcal{S}}
\def\T{\mathcal{T}}
\def\P{\mathcal{P}}
\def\I{\mathcal{I}}
\def\J{\mathcal{J}}
\def\Q{\mathcal{Q}}
\def\U{\mathcal{U}}

\def\g{\mathfrak g}
\def\z{\mathfrak z}
\def\m{\mathfrak m}

\def\vol{{\rm vol\ }}
\def\Im{{\rm Im\ }}
\def\Ker{{\rm Ker\ }}
\def\hol{{\rm hol}}
\def\Aut{{\rm Aut}}
\def\End{{\rm End}}
\def\Pr{{\rm Pr}}
\def\pr{{\rm pr}}
\def\id{{\rm id}}
\def\ad{{\rm ad}}
\def\Ad{{\rm Ad}}
\def\rel{{\rm rel}}
\def\ord{{\rm ord}}

\def\bv{{\bf v} }
\def\bV{{\bf V} }
\def\ba{{\bf a} }
\def\bb{{\bf b} }
\def\bc{{\bf c} }
\def\bu{{\bf u} }
\def\bw{{\bf w} }
\def\bx{{\bf x} }
\def\bW{{\bf W} }
\def\bk{{\bf k} }
\def\b0{{\bf 0} }

\def\eps{\varepsilon}

\begin{document}

\vskip 1.0 true cm

\begin{center}
{\Huge \bf Stability of the space of }

\vskip 0.4 true cm

{\Huge \bf Automorphic Forms under }

\vskip 0.4 true cm

{\Huge \bf Local Deformations of the }

\vskip 0.4 true cm

{\Huge \bf Lattice}
\end{center}

\vskip 1.0 true cm

\begin{center}
R. Knevel

{\it Université du Luxembourg, Unité de Recherche en Mathématiques}

({\it Received 28 May 2009})
\end{center}

\begin{small}


\begin{quote}


{\bf Abstract:} First we explain the concept of local deformation over a 'parameter' algebra $\P$ , in particular the notion of a $\P$-lattice in a {\sc Lie} group. The purpose of this article is to define the spaces $M_k(\Upsilon)$ and 
$S_k(\Upsilon)$ of automorphic resp. cusp forms on the upper half plane $H$ for a $\P$- (!) lattice $\Upsilon$ of $SL(2, \rz)$ and to investigate their structure. It turns out that in almost all cases the spaces $M_k(\Upsilon)$ and $S_k(\Upsilon)$ are 
free modules over the complexified $\P$ of rank equal to the dimension of the spaces of automorphic resp. cusp forms for the body $\Gamma := \Upsilon^\#$ , which is the associated ordinary lattice in $SL(2, \rz)$~. In other words almost every 
automorphic resp. cusp form admits an 'adaption' to local deformations of the lattice. This is shown by giving the quotient $\Upsilon \backslash H$ together with the cusps of $\Gamma \backslash H$ the structure of a $\P$- {\sc Riemann} surface and 
writing the spaces of automorphic resp. cusp forms as global sections of holomorphic $\P$- (!) line bundles on $\Upsilon \backslash H \cup \{ \text{ cusps of } \Gamma \backslash H \}$~.

\end{quote}

\end{small}

\section{Introduction}

First of all let us discuss the concept of local deformation. A rough explanation is the following: Vary the 'data' describing a classical object, for example the glueing data of local charts defining a 
smooth manifold, let them depend on 'parameters' generating a local commutative algebra $\P$ . In practice it is not necessary to specify the parameters since all information is already encoded in this algebra. Here we consider the case 
of an finite dimensional algebra $\P$ whose unique maximal ideal $\I$ is nilpotent, and $\P / \I \simeq \rz$ . This case obviously lies in-between the extremal poles of infinitesimal deformation, which means $\I^2 = 0$ , and arbitrary local deformation. So it 
is not surprising that many proofs in this context use the techniques coming from infinitesimal deformation in combination with induction over the power annihilating $\I$ . An object $\O$ whose data depend on the 'parameters' generating $\P$ will be called 
a $\P$-object. It is necessarily a local deformation of a classical object $\O^\#$ , called its body. We will see that in general this gives whole body functors from the categories of $\P$-objects to their classical counterparts induced by the canonical projection 
$\P \rightarrow \P / \I \simeq \rz$ , heuristically the 'set all the parameters to $0$ ' functors. \\

Natural questions in the framework of local deformation are the following: Do we really get more objects when we allow $\P$-ones, in other words do certain objects allow non-trivial local deformation or are they completely rigid? Can we adapt 'functions' 
on a classical object $O$ to local deformations of $O$ ? How can one classify all $\P$-objects with given body? \\

In this article we study $\P$-lattices in $SL(2, \rz)$ acting on the upper half plane $H \subset \cz$ , in other words local deformation of the natural embedding of a given lattice $\Gamma \hookrightarrow SL(2, \rz)$ as group homomorphism and we want to 
investigate the spaces of automorphic and cusp forms for these $\P$-lattices. The theory of automorphic forms for classical lattices is already well-established. Let $\Gamma \sqsubset SL(2, \rz)$ be a lattice. Then we have an asymptotic formula

\begin{equation*}
\dim M_k(\Gamma) , \dim S_k(\Gamma) \sim \frac{k}{4 \pi} \vol(\Gamma \backslash H)
\end{equation*} \\

for the dimension of the spaces $M_k(\Gamma)$ and $S_k(\Gamma)$ of automorphic resp. cusp forms for $\Gamma$ of high weight $k$ , and this is one of the most beautiful applications of the theory of holomorphic line bundles on compact 
{\sc Riemann} surfaces. Now in the case of a $\P$-lattice $\Upsilon$ in $SL(2, \rz)$ with body $\Gamma$ it would be nice if every form $f \in M_k(\Gamma)$ would allow an adaption $\widetilde f \in M_k(\Upsilon)$ to $\Upsilon$ having $f$ as 
body, because this is equivalent to the stability of $M_k(\Gamma)$ under local deformations of $\Gamma$ . We will show that this is precisely equivalent to $M_k(\Upsilon) \simeq M_k(\Gamma) \otimes \P^\cz$ as $\P^\cz$-modules and, as the main result 
of this paper, is indeed true except in the case where $k = 1$ and $\Gamma \backslash H \cup \{ \text{ cusps of } \Gamma \backslash H \}$ is of genus $g \geq 2$ , see theorems \ref{main} and \ref{main2} for even resp. odd $k \in \nz$ . \\

Imitating the classical theory, the article is organized as follows: The general concept of $\P$-manifolds and $\P$- vector bundles is introduced in section \ref{param} , and section \ref{lat and aut} treats the basic definitions of $\P$-lattices and associated automorphic 
and cusp forms in the case of $SL(2, \rz)$ acting on $H$~. In section \ref{quot} after fixing a $\P$-lattice $\Upsilon$ of $SL(2, \rz)$ with body $\Gamma := \Upsilon^\#$ we construct a $\P$- {\sc Riemann} 
surface $\X$ which is in some sense a compactification of $\Upsilon \backslash H$ in analogy to $X := \Gamma \backslash H \cup \{ \text{ cusps of } \Gamma \backslash H\}$ , which will be the body of $\X$~. In section \ref{Teich} we do some elementary 
{\sc Teichmüller } theory, more precisely we prove that any $\P$- {\sc Riemann} surface $\X$ with compact body $X := \X^\#$ is represented by a $\P$-point of the {\sc Teichmüller} space whose body represents $X$ , see theorem 
\ref{paramTeich} . This result is of course of general interest since it gives a complete solution for the classification problem of $\P$- {\sc Riemann} surfaces with given compact body and is therefore given in greatest possible generality. In section \ref{sect main}
we use all our knowledge obtained so far to prove the main theorems of this article, theorems \ref{main} and \ref{main2} . Finally, section \ref{SL} deals with the special case $\Upsilon^\# = SL(2, \zz)$ . \\

\section{$\mathcal{P}$-manifolds} \label{param}

For the whole article let $\P$ be a finite dimensional real unital commutative algebra with a unital algebra projection ${}^\# : \P \rightarrow \rz$ and the unique maximal ideal $\I := \Ker {}^\# \lhd \P$ such that $\I^N = 0$ for 
some $N \in \nz$ . ${}^\#$ is called the body map of $\P$ .

\begin{examples}
\end{examples}

\begin{itemize}
\item[(i)] The even part $\P := \bigwedge\rund{\rz^{N - 1}}_0$ of an exterior algebra,
\item[(ii)] the polynomial algebra $\P := \rz[X] \left/ \rund{X^N = 0}\right.$ with cut off.
\end{itemize}

Defining the category of $\P$-manifolds will be done in terms of ringed spaces, the real and complex case treated simultaneously. Therefore let $\P^\cz$ and $\I^\cz$ denote the complexifications of $\P$ resp. $\I$ .

\begin{defin} [ $\P$-manifolds and $\P$-points]
\item[(i)] Let $M$ be a real smooth (complex) manifold of dimension $n$ , and $\S$ be a sheaf of unital commutative $\P$- ( $\P^\cz$- ) algebras over $M$ . Then the ringed space $\M := (M, \S)$ is called a real (complex) $\P$-manifold of dimension $n$ if 
and only if locally $\S \simeq \C^\infty_M \otimes \P$ ( $\S \simeq \O_M \otimes \P^\cz$ ) . $\M^\# := M$ is called the body of $\M$ . If $M$ is a complex manifold of dimension $n = 1$ then $\M$ is called a $\P$- {\sc Riemann} surface.

\item[(ii)] Let $\M = (M, \S)$ and $\N = (N, \T)$ be two real (complex) $\P$-manifolds. A $\P$-morphism between $\M$ and $\N$ is a morphism from $\M$ two $\N$ as ringed spaces, more precisely a collection 
$\Phi := \rund{\varphi, \rund{\phi_V}_{V \subset N \text{ open } }}$ where \\
$\varphi: M \rightarrow N$ is a smooth (holomorphic) map, and all \\
$\phi_V : \T(V) \rightarrow \S\rund{\varphi^{- 1}(V)}$ are unital $\P$- ( $\P^\cz$- ) algebra homomorphisms such that for all $W \subset V \subset N$ open

\[
\begin{array}{ccc}
\phantom{123} \T(V) & \mathop{\longrightarrow}\limits^{\phi_V} & \S\rund{\varphi^{- 1}(V)} \phantom{12345} \\
|_W \downarrow & \circlearrowleft & \downarrow |_{\varphi^{- 1}(W)} \\
\phantom{123} \T(W) & \mathop{\longrightarrow}\limits_{\phi_W} & \S\rund{\varphi^{- 1}(W)} \phantom{12345}
\end{array}  \, .
\]

We write $\Phi: \M \rightarrow_{\P} \N$ and $f(\Phi) := \phi_V(f) \in \S\rund{\varphi^{- 1}(V)}$ for all $V \subset N$ open and $f \in \T(V)$ . $\Phi^\# := \varphi$ is called the body of $\Phi$ .

\item[(iii)] Let $a = \rund{a^\#, \rund{a_V}_{V \subset N \text{ open } }} : \M \rightarrow_{\P} \N$ be a $\P$-morphism from the real (complex) $\P$-manifold $\M := \rund{\{0\}, \P}$ ( $\M := \rund{\{0\}, \P^\cz}$ ) to the real (complex) $\P$-manifold $\N$ . Then $a$ 
is called a $\P$-point of $\N$ . Its body $a^\# : \{0\} \hookrightarrow \N^\#$ will always be identified with the usual point $a^\#(0) \in \N^\#$~. We write 
$a \in_{\P} \N$ and $f(a) := a_V(f)$ for all $V \subset N$ open with $a^\# \in V$ and $f \in \T(V)$ . If $\O$ is another $\P$-manifold and $\Phi: \N \rightarrow_{\P} \O$ a $\P$-morphism then we write $\Phi(a) := \Phi \circ a \in_{\P} \O$ .
\end{defin}

Let us collect some basic properties of $\P$-manifolds, they will be implicitly used later:

\begin{itemize}
\item[(i)] Let $\M = (M, \S)$ be a real (complex) $\P$-manifold. Then since $\C^\infty_M$ (~$\O_M$ ) has no other unital sheaf automorphisms than $\id$ , ${}^\# : \P \rightarrow \rz$ induces a body map ${}^\# : \S \rightarrow \C^\infty_M$ 
( ${}^\# : \S \rightarrow \O_M$ ) , which is a projection of sheaves of real (complex) unital algebras. The kernel of ${}^\#$ is precisely $\I \S$ ( $\I^\cz \S$ ), and we have a canonical sheaf isomorphism 
$\C^\infty_M \otimes \I^n \simeq \I^n \S$ ( $\O_M \otimes \rund{\I^\cz}^n \simeq \rund{\I^\cz}^n \S$ ) whenever $n \in \nz$ such that $\I^{n + 1} = 0$ . \\

Now let $\N := (N, \T)$ be another real (complex) $\P$-manifold and \\
$\Phi := \rund{\varphi, \rund{\phi_V}_{V \subset N \text{ open } }}$ a $\P$ morphism from $\M$ to $\N$ . Then automatically

\[
\begin{array}{ccc}
\phantom{12} \T(V) & \mathop{\longrightarrow}\limits^{\phi_V} & \S\rund{\varphi^{- 1}(V)} \phantom{1} \\
{}^\# \downarrow & \circlearrowleft & \downarrow {}^\# \\
\phantom{12} \begin{array}{c} \C^\infty(V) \\
( \phantom{12} \O(V) \phantom{12} ) \end{array} & \mathop{\longrightarrow}\limits_{f \mapsto f \circ \varphi|_{\varphi^{- 1}(V)} } & \begin{array}{c} \C^\infty\rund{\varphi^{- 1}(V)} \\
( \phantom{12} \O\rund{\varphi^{- 1}(V)} \phantom{12} ) \end{array}
\end{array}
\]

for all $V \subset N$ open, and for all $f \in \T(V)$ we call \\
$\phi_V(f) = f(\Phi) \in \S\rund{\varphi^{- 1}(V)}$ the pullback of $f$ under $\Phi$ .

\item[(ii)] Every usual real smooth (complex) manifold $M$ can be regarded as a real (complex) $\P$-manifold identifying $M$ with the ringed space $\rund{M, \C^\infty_M \otimes \P}$ ( $\rund{M, \O_M \otimes \P^\cz}$ ) , and every usual smooth 
(holomorphic) map between real smooth (complex) manifolds can be regarded as a $\P$-morphism between them.

\item[(iii)] Let $\Phi = \rund{\varphi, \rund{\phi_V}_{V \subset N \text{ open } }} : \M \rightarrow_{\P} \N$ be a $\P$-morphism between the real (complex) $\P$-manifolds $\M = (M, \S)$ and $\N = (N, \T)$ . Then it is an isomorphism iff $\varphi$ is 
bijective, and in this case $\Phi$ is called a $\P$-isomorphism. \\

If $\varphi$ is an immersion then for all $a \in M$ there exists an open neighbourhood $V$ of $\varphi(a)$ in $N$ such that $\phi_V : \T(V) \rightarrow \S\rund{\varphi^{- 1}(V)}$ is surjective. $\Phi$ is called a $\P$-embedding iff $\varphi$ is an 
injective immersion and so a smooth (holomorphic) embedding of real smooth (complex) manifolds, and in this case we write $\Phi: \M \hookrightarrow_{\P} \N$ . \\

Furthermore, if $\varphi$ is surjective and for all $a \in M$ there exists an open neighbourhood $V$ of $\varphi(a)$ in $N$ such that $\phi_V : \T(V) \rightarrow \S\rund{\varphi^{- 1}(V)}$ is injective then $\Phi$ is called a $\P$-projection.

\item[(iv)] Let $\Phi = \rund{\varphi, \rund{\phi_V}_{V \subset N \text{ open } }} : \M \hookrightarrow_{\P} \N$ be a $\P$-embedding from the real (complex) $\P$-manifold $\M = (M, \S)$ into the real 
(complex) $\P$-manifold $\N = (N, \T)$ . Since $\varphi: M \hookrightarrow N$ is an embedding we can regard $M$ as a real smooth (complex) submanifold of $N$ and $\varphi$ as the 
canonical inclusion. So $\Phi$ induces a sheaf projection $\pi: \T|_M \rightarrow \S$ such that for all $V \subset N$ open

\[
\begin{array}{ccc}
\T(V) & \mathop{\longrightarrow}\limits^{\phi_V} & \S(V \cap M) \\
{}_{|_{V \cap M}} \searrow & \circlearrowleft & \nearrow_{\pi_{V \cap M}} \\
 & \T|_N(V \cap M) & 
\end{array} \, .
\]

In particular there is a 1-1 correspondence between $\P$-points \\
$a = \rund{a^\#, \rund{a_V}_{V \subset N \text{ open } }}$ of $\N$ and pairs $\rund{a^\#, a_{a^\#}}$ where $a^\# \in \N^\#$ and $a_{a^\#} : \T_{a^\#} \rightarrow \P$ ( $a_{a^\#} : \T_{a^\#} \rightarrow \P^\cz$ ) is an epimorphism of $\P$- (~$\P^\cz$- ) algebras, 
where $\T_{a^\#}$ denotes the stalk of $\T$ at $a^\# \in N$ , such that for all $V \subset N$ open if $a^\# \in V$ then $a_V = a_{a^\#} \circ |_{a^\#}$ , where $|_{a^\#} : \T(V) \rightarrow \T_{a^\#}$ denotes the canonical projection, and otherwise $a_V \equiv 0$ .

\item[(v)] The local models of real (complex) $\P$-manifolds are the usual open sets $U \subset \rz^m$ ( $U \subset \cz^m$ ) regarded as real (complex) $\P$-manifolds together with $\P$-morphisms between them. Let $\M = (M, \S)$ be an $m$-dimensional 
real (complex) $\P$-manifold and $V \subset \rz^n$ ( $V \subset \cz^n$ ) be open. Then one can show that there is a 1-1 correspondence between $\P$-morphisms from $U$ to $V$ and $n$-tuples $\rund{f_1, \dots, f_n} \in \S(M)^{\oplus n}$ such that 
$\rund{f_1^\#(u), \dots, f_n^\#(u)} \in V$ for all $u \in M$ given as follows:

\begin{quote}
To a $\P$-morphism $\Phi := \rund{\varphi, \rund{\phi_W}_{W \subset V \text{ open } }} : \M \rightarrow_{\P} V$ we associate the tuple $\rund{\phi_V\rund{x_1}, \dots, \phi_V\rund{x_n}}$ \\
( $\rund{\phi_V\rund{z_1}, \dots, \phi_V\rund{z_n}}$ ) , where $x_1, \dots, x_n \in \C^\infty(V)$ \\
( $z_1, \dots, z_n \in \O(V)$ ) denote the coordinate functions on $V$ . \\

Conversely to a tuple $\rund{f_1, \dots, f_n}$ we associate the $\P$-morphism $\Phi := \rund{\varphi, \rund{\phi_W}_{W \subset V \text{ open } }} : \M \rightarrow_{\P} V$ , where $\varphi: M \rightarrow V \, , \, u \mapsto \rund{f_1^\#(u), \dots, f_n^\#(u)}$ and

\begin{eqnarray*}
&& \phi_W : \C^\infty(W) \otimes \P \rightarrow \S\rund{\varphi^{- 1}(W)} \\
&& \phantom{1234567} ( \phantom{12} \O(W) \otimes \P^\cz \rightarrow \S\rund{\varphi^{- 1}(W)} \phantom{12} ) \, , \\
&& \phantom{12} h \mapsto \sum_{\bk \in \nz^n} \frac{1}{\bk !} \rund{\rund{\partial^{\bk} h} \circ \rund{\varphi|_{\varphi^{- 1}(W)}} } \times \\
&& \phantom{1234567} \times \left.\rund{f_1 - f_1^\#}^{k_1} \cdots \rund{f_n - f_n^\#}^{k_1} \right|_{\varphi^{- 1}(W)}
\end{eqnarray*}

for all $W \subset V$ open, which is nothing but the formal {\sc Taylor} expansion of the expression $h\rund{f_1, \dots, f_n}$ .
\end{quote}

In particular one can identify the $\P$-points $\ba$ of $V$ with the tuples $\rund{a_1, \dots, a_n} \in \P^{\oplus n}$ ( $\rund{a_1, \dots, a_n} \in \rund{\P^\cz}^{\oplus n}$ ) such that $\rund{a_1^\#, \dots, a_n^\#} \in V$ as follows:

\begin{quote}
To a $\P$-point $\ba \in_{\P} V$ one assigns the tuple $\rund{\ba_V\rund{x_1}, \dots, \ba_V\rund{x_n}}$ (resp. $\rund{\ba_V\rund{z_1}, \dots, \ba_V\rund{z_n}}$ ), \\

and conversely to a tuple $\rund{a_1, \dots a_n}$ one assigns \\
$\ba = \rund{\ba^\#, \ba_{\ba^\#}}$ , where $\ba^\# = \rund{a_1^\#, \dots a_n^\#} \in V$ and

\begin{eqnarray*}
&& \ba_{\ba^\#} : \left.\C^\infty_V \right|_{\ba^\#} \otimes \P \rightarrow \P \phantom{12345} ( \phantom{12} \left.\O_V \right|_{\ba^\#} \otimes \P^\cz \rightarrow \P^\cz \phantom{12} ) \, , \\
&& \phantom{12} h \mapsto \sum_{\bk \in \nz^n} \frac{1}{\bk !} \partial^{\bk} h \rund{\ba^\#} \rund{a_1 - a_1^\#}^{k_1} \cdots \rund{a_n - a_n^\#}^{k_1}  \, ,
\end{eqnarray*}

which is again nothing but the formal {\sc Taylor} expansion of the expression $h\rund{a_1, \dots, a_n}$ .
\end{quote}

Now let $\M = (M, \S)$ and $\N = (N, \T)$ be real (complex) $\P$-manifolds and $V \subset \rz^k$ ( $V \subset \cz^k$ ) open. Let $\Phi: \M \rightarrow_{\P} \N$ and $\Psi: \N \rightarrow_{\P} V$ be $\P$-morphisms and $\Psi$ be given by the tuple 
$\rund{f_1, \dots, f_k} \in \T(N)^{\oplus k}$~. Then $\Psi \circ \Phi$ is given by the tuple $\rund{f_1(\Phi), \dots, f_k(\Phi)} \in \S(M)^{\oplus k}$~. In particular if $a \in_{\P} \N$ then $\Psi(a) \in_{\P} V$ is given by the tuple 
$\rund{f_1(a), \dots, f_k(a)} \in \P^{\oplus k}$ ( $\rund{f_1(a), \dots, f_k(a)} \in \rund{\P^\cz}^{\oplus k}$ ). \\

Of course given a real (complex) $\P$-manifold $\M = (M, \S)$ of dimension $n$ for each $a \in M$ there exists an open neighbourhood $U \subset M$ of $a$ , $V \subset \rz^n$ ( $V \subset \cz^n$ ) open and a $\P$-isomorphism from $(U, \S|_U)$ to 
$\rund{V, \C^\infty_V \otimes \P}$ ( $\rund{V, \O_V \otimes \P}$ ) . Such a $\P$-isomorphism is called a local $\P$-chart of $\M$ at $a$ . Two local $\P$-charts $V_i , V_j \subset \rz^n$ (~$V_i , V_j \subset \cz^n$~) 'glue' together via a $\P$-glueing data, given as 
a $\P$-isomorphism $\Phi_{i j} : V_{i j} \rightarrow_{\P} V_{j i}$ between the overlaps $V_{i j} \subset V_i$ and $V_{j i} \subset V_j$ open. The body then will be given by the same local charts with glueing data $\Phi_{i j}^\#$ . Observe that in general one can 
not specify ordinary points of a real (complex) $\P$-manifold, only $\P$-points. However, given a $\P$-point $a \in_{\P} \M$ there always exists a local $\P$-chart of $\M$ at $a^\#$ mapping $a$ to a usual point of $\rz^n$ ( $\cz^n$ ). \\

Now let $U \subset \cz^m$ and $V \subset \cz^n$ be open. Then $U$ and $V$ can be regarded as open subsets of $\rz^{2 m}$ resp. $\rz^{2 n}$ , and using the 1-1 correspondence from above we see that every $\P$-morphism from $U$ to $V$ 
regarded as complex open sets is at the same time a $\P$-morphism from $U$ to $V$ as real open sets, and so we get a whole 'forget' functor from the category of complex $\P$-manifolds to the one of real $\P$-manifolds. \\

\item[(vi)] Let $\M$ be a real $\P$-manifold. Then there exists a $\P$-isomorphism $\Phi: \M^\# \rightarrow_{\P} \M$ such that $\Phi^\# = \id$ . This can be shown by induction on $N$ using $H^1\rund{\M^\#, T \M^\#} = 0$ , and it is nothing but the rigidity of 
smooth manifolds under local deformation. \\

\item[(vii)] Let $\M = (M, \S)$ be a real (complex) $\P$-manifold of dimension $m$~, \\
$\Phi := \rund{\varphi, \rund{\Phi_V}_{V \subset M \text{ open } }} : \M \rightarrow_{\P} \rz^n$ ( $\M \rightarrow_{\P} \cz^n$ ) be a $\P$-morphism such that $D \Phi^\#$ is 
surjective at every point of $M$ and Finally, $\ba \in_{\P} \rz^n$ (~$\ba \in_{\P} \cz^n$ ) . Then we can define the preimage $\Phi^{- 1}(\ba)$ of $\ba$ under $\Phi$ as a real (complex) $\P$-manifold $(N, \T)$ of dimension $m - n$ as follows: Its body is
$N := \varphi^{- 1}\rund{\ba^\#}$ , and the sheaf $\T$ is given by

\[
\T := \left.\S|_N \, \right/ \m \, ,
\]

where $\m \lhd \S|_N$ is the ideal sheaf generated by all $\phi_V(f)$ , where $V \subset \rz^n$ (~$V \subset \cz^n$ ) open and $f \in \C^\infty(V)$ ( $f \in \O(V)$ ) such that $f(\ba) = 0$ if $\ba^\# \in V$ . If $\Phi$ is given by the tuple 
$\rund{f_1, \dots, f_n} \in \S(M)^{\oplus n}$ and $\ba$ by $\rund{a_1, \dots, a_n} \in \P^{\oplus n}$ ( $\rund{a_1, \dots, a_n} \in \rund{\P^\cz}^{\oplus n}$ ) then $\m$ is generated by $f_1 - a_1 , \dots , f_n - a_n$ . The $\P$-morphism 
$I := \rund{i, \rund{\pi_U}_{U \subset M \text{ open } }} : \N \hookrightarrow_{\P} \M$ , where $i: N \hookrightarrow M$ is the canonical inclusion, and $\pi_{U} : \S(U) \rightarrow \T(U \cap N)$ , $U \subset M$ open, denote the canonical projections, is a 
$\P$-embedding called the canonical inclusion of $\N$ into $\M$ . \\

Let $\O$ be another real (complex) $\P$-manifold. Then there exists a 1-1 correspondence between the $\P$-morphisms $\Psi: \O \rightarrow_{\P} \N$ and the $\P$-morphisms $\Xi: \O \rightarrow_{\P} \M$ having

\[
\begin{array}{ccc}
\phantom{123} \O & \mathop{\longrightarrow}\limits^{\Xi} & \M \phantom{12} \\
\Pr \downarrow & \circlearrowleft & \downarrow \Phi \\
\phantom{123} \{0\} & \mathop{\longrightarrow}\limits_{\ba} & \rz^n \phantom{1}
\end{array} \, ,
\]

where $\Pr: \O \rightarrow_{\P} \{0\}$ denotes the canonical $\P$-projection. It is given by the assignment $\Psi \mapsto I \circ \Psi$ . \\

In particular we can identify the $\P$-points of $\N$ with the $\P$-points $b \in_{\P} \M$ of $\M$ having $\Phi(b) = \ba$ .

\item[(viii)] In the category of real (complex) $\P$-manifolds there exists a cross product. If $\M = (M, \S)$ and $\N = (N, \T)$ are two real (complex) $\P$-manifolds then their cross product is given by

\[
\M \times \N := \rund{M \times N, \pr_1^\ast \S \hat\otimes \pr_2^\ast \T} \, ,
\]

and the canonical $\P$-projections by \\
$\Pr_1 := \rund{\pr_1, \rund{i_U}_{U \subset M \text{ open } }} : \M \times \N \rightarrow_{\P} \M$ and \\
$\Pr_2 := \rund{\pr_2, \rund{j_V}_{V \subset N \text{ open } }} : \M \times \N \rightarrow_{\P} \N$ , where

\[
i_U : \S(U) \hookrightarrow \rund{\pr_1^\ast \S \hat\otimes \pr_2^\ast \T}\rund{\pr_1^{- 1}(U)} = \S(U) \hat\otimes \T(N)
\]

and

\[
j_V : \T(V) \hookrightarrow \rund{\pr_1^\ast \S \hat\otimes \pr_2^\ast \T}\rund{\pr_2^{- 1}(V)} = \S(M) \hat\otimes \T(V) \, ,
\]

$U \subset M$ , $V \subset N$ open, denote the canonical inclusions. \\

By the universal property of the cross product there is a 1-1 correspondence between the $\P$-points $c \in_{\P} \M \times \N$ of $\M \times \N$ and pairs $(a, b)$ of $\P$-points $a \in_{P} \M$ and $b \in_{\P} \N$ given by the assignment
$c \mapsto \rund{\Pr_1(c), \Pr_2(c)}$ .
\end{itemize}

Let me give two typical proofs: \\

First we prove the statement of (iii) in the real case (same proof in the complex case): "$\Rightarrow$" is trivial. "$\Leftarrow$" will be proven by induction on $n \in \nz$ with $\I^n = 0$ , $\I \lhd \P$ being the unique maximal ideal of $\P$ . If $n = 1$ then of 
course the statement is trivial since then $\Phi_V(f) = f \circ \varphi|_{\varphi^{- 1}(V)}$ for all $V \subset N$ open and $f \in \C^\infty(V)$ . \\

Now let $\I^{n + 1} = 0$ . Then we define $\Q := \P \left/ \I^n\right.$ . Clearly $\Q$ has \\
$\J := \I \left/ \I^n\right.$ as unique maximal ideal, and $\J^n = 0$ . Let ${}^\natural : \P \rightarrow \Q$ denote the canonical projection. $\M^\natural := \rund{M, \S \left/ \I^n \S\right.}$ and $\N^\natural := \rund{N, \T \left/ \I^n \T\right.}$ are real $\Q$-manifolds, 
and ${}^\natural$ induces sheaf projections

\[
{}^\natural: \S \rightarrow \S \left/ \I^n \S\right.
\]

and

\[
{}^\natural: \T \rightarrow \T \left/ \I^n \T\right.  \, .
\]

Now let $\Phi = \rund{\varphi, \rund{\phi_V}_{V \subset N \text{ open }} }$ be a $\P$-morphism from $\M$ to $\N$ . Then $\Phi$ induces a $\Q$-morphism $\Phi = \rund{\varphi, \rund{\phi_V^\natural}_{V \subset N \text{ open }} }$ from $\M^\natural$ to 
$\N^\natural$~, where for all $V \subset N$ open $\phi_V^\natural$ is the unique unital $\P$-algebra morphism $\rund{\T \left/ \I^n \T\right.}(V) \rightarrow \rund{\S \left/ \I^n \S\right.}\rund{\varphi^{- 1}(V)}$ such that

\[
\begin{array}{ccc}
\phantom{12} \T(V) & \mathop{\longrightarrow}\limits^{\phi_V} & \S\rund{\varphi^{- 1}(V)} \\
{}^\natural \downarrow & \circlearrowleft & \downarrow {}^\natural \\
\phantom{12} \rund{\T \left/ \I^n \T\right.}(V) & \mathop{\longrightarrow}\limits_{\phi_V^\natural} & \rund{\S \left/ \I^n \S\right.}\rund{\varphi^{- 1}(V)}
\end{array}  \, .
\]

Now we have to show that $\Phi$ is an isomorphism of ringed spaces. But since $\varphi$ is already bijective it suffices to show that $\Phi$ is a local isomorphism. Therefore we may assume without loss of generality that $\S = \C^\infty_M \otimes \P$ 
and $\T = \C^\infty_N \otimes \P$ . So we know that $\Phi_V(h) = h \circ \varphi|_{\varphi^{- 1}(V)}$ for all $V \subset N$ open and $h \in \C^\infty(V) \otimes \I^n \hookrightarrow \T(V)$ . \\

By induction hypothesis we already know that $\Phi^\natural$ is a $\P$-isomorphism. Let $V \subset N$ be open and $h \in \S\rund{\varphi^{- 1}(V)}$ . Then $\phi_V^\natural$ is an isomorphism, and so there exists $f  \in \T(V)$ such that 
$\phi_V^\natural\rund{f^\natural} = h^\natural$ , and therefore \\
$\Delta := h - \phi_V(f) \in \C^\infty\rund{\varphi^{- 1}(V)} \otimes \I^n$ . Since $\varphi: M \rightarrow N$ is a diffeomorphism we can build $\Delta \circ \left.\varphi^{- 1}\right|_V \in \C^\infty(V) \otimes \I^n$ , and

\[
\phi_V\rund{f + \Delta \circ \left.\varphi^{- 1}\right|_V} = h - \Delta + \phi_V\rund{\Delta \circ \left.\varphi^{- 1}\right|_V} = h
\]

by (i) . This proves surjectivity of $\phi_V$ . For proving injectivity let $f \in \T(V)$ such that $\Phi_V(f) = 0$ . Then $\Phi_V^\natural\rund{f^\natural} = 0$ , and so $f^\natural = 0$ . Therefore $f \in \C^\infty(V) \otimes \I^n$ , and so
$0 = \Phi_V(f) = f \circ \varphi|_{\varphi^{- 1}(V)}$ . This implies $f = 0$ . $\Box$ \\

Now we prove that $(N, \T)$ in (vii) is indeed a $\P$-manifold of dimension $m - n$ in the real case (it is again the same proof in the complex case): Let $\bx_0 \in N$~. Then it is enough to show that there exists an open neighbourhood $U \subset M$ of 
$\bx_0$ such that $\rund{U \cap N, \T|_{U \cap N}}$ is a real $\P$-manifold. So first of all choose a neighbourhood $U \subset M$ of $\bx_0$ such that $\rund{U \cap N, \T|_{U \cap N}}$ is identified with an open subset of $\rz^m$ regarded as real 
$\P$-manifold, and without loss of generality we may assume that $M = U$ . So let $\Phi$ be given by the tuple $\rund{f_1, \dots, f_n} \in \rund{\C^\infty(U) \otimes \P}^{\oplus n}$ and $\ba$ by the tuple $\rund{a_1, \dots, a_n} \in \rund{\P^\cz}^{\oplus n}$ . 
Then after maybe replacing $\rund{f_1, \dots, f_n}$ by $\rund{f_1 - a_1, \dots, f_n - a_n}$ we may assume without loss of generality that $\ba = \b0$ . Now let the $\P$-morphism $\widetilde \Phi$ from $U$ to $\rz^m$ be given by the tuple 
$\rund{f_1, \dots, f_n, x_{n + 1}, \dots, x_m}$ , where $x_1, \dots, x_m \in \C^\infty(U)$ denote the standard coordinate functions on $U$ . Then since $D \Phi^\# \rund{\bx}$ is surjective at every point $\bx \in U$ by assumption after maybe changing the 
order of the coordinates we assume without loss of generality that $D {\widetilde \Phi}^\# \rund{\bx_0} \in GL(m, \rz)$ . So after maybe replacing $U$ by a smaller open neighbourhood of $\bx_0$ we may assume without loss of generality 
that $\varphi := {\widetilde \Phi}^\#$ is a diffeomorphism from $U$ to $V := \varphi(U) \subset \rz^m$ , and so by (iii) $\widetilde \Phi$ is a $\P$-isomorphism from $U$ to $V$ . But then we see that $\Phi \circ {\widetilde \Phi}^{- 1}$ is precisely the projection 
onto the first $n$ coordinates, which is a usual smooth map from $V$ to $\rz^n$ . So identifying $U$ and $V$ via $\widetilde \Phi$ we may without loss of generality assume that $\widetilde \Phi = \id$ , and then the statement is trivial. $\Box$ \\

Let us already here introduce the notion of $\P$- vector bundles over $\P$-manifolds. It will be crucial in section \ref{sect main} .

\begin{defin} [ $\P$- vector bundles]
\item[(i)] Let $\M = (M, \S)$ be a real (complex) $\P$-manifold. Then an $\S$- sheaf module $\E$ on $M$ is called a real (holomorphic) $\P$- vector bundle of rank $r$ over $\M$ iff it is locally isomorphic to $\S^{\oplus r}$ . In this case 
${}^\# : \P \rightarrow \rz$ induces a body map ${}^\# : \E \rightarrow \Gamma^{\infty}(\diamondsuit, E)$ ( ${}^\# : \E \rightarrow \Gamma^{\hol}(\diamondsuit, E)$ ), where $\Gamma^{\infty}(\diamondsuit, E)$ 
( ${}^\# : \E \rightarrow \Gamma^{\hol}(\diamondsuit, E)$ ) is the sheaf of smooth (holomorphic) sections of a real smooth (holomorphic) vector bundle $E \rightarrow M$ is of rank $r$ , which is uniquely determined by $\E$ . $\E^\# := E$ is called the body of 
$\E$ . The space $H^0(\E) := \E(M)$ is called the space of global sections of $\E$ , it is a $\P$- ( $\P^\cz$-~) module. If $r = 1$ then $\E$ is called a $\P$- line bundle.

\item[(ii)] Let $\E$ be a real (holomorphic) $\P$- vector bundle over the real (complex) $\P$-manifold $\M = (M, \S)$ , and let $\Phi: \N \hookrightarrow \M$ be a $\P$-embedding of the real (complex) $\P$-manifold $\N = (N, \T)$ .

\[
\E|_{\N} := \left.\E|_N \frac{}{} \right/ \m \, \E|_N \, ,
\]

where $\m$ denotes the kernel of the canonical sheaf projection $\S|_N \rightarrow \T$ , is called the restriction of the $\P$- vector bundle $\E$ to $\N$ . It is a real (holomorphic) $\P$- vector bundle over $\N$ of rank $r$ with body $E|_N$ . If $U \subset M$ 
open and $F \in \E(U)$ then the image $F|_{\N}$ of $F$ under the canonical map $\E(U) \rightarrow \E|_N\rund{\varphi^{- 1}(U)} \rightarrow \E|_{\N}\rund{\varphi^{- 1}(U)}$ is called the restriction of $F$ to $\N$ .

\item[(iii)] Let $\E$ and $\F$ be real (holomorphic) $\P$- vector bundles over the real (complex) $\P$-manifold $\M = (M, \S)$ of rank $r$ resp. $s$ with bodies $E$ resp. $F$ . Then $\E \otimes \F := \E \otimes_{\S} \F$ is called the tensor product of $\E$ and 
$\F$ . It is a $\P$- vector bundle over $\M$ of rank $r s$ with body $E \otimes F$ .
\end{defin}

Let us collect some basic facts about $\P$- vector bundles:

\begin{itemize}
\item[(i)] If $E$ is a usual real smooth (holomorphic) vector bundle of rank $r$ over the usual real smooth (complex) manifold $M$ then $E$ can be identified with the real (holomorphic) $\P$- vector bundle 
$\Gamma^\infty(\diamondsuit, E) \otimes \P$ ( $\Gamma^\hol(\diamondsuit, E) \otimes \P^\cz$ ) over $M$ regarded as the real (complex) $\P$-manifold $(M, \C^\infty_M \otimes \P)$ ( $(M, \O_M \otimes \P^\cz)$ ) . \\

\item[(ii)] Let $\E$ be the real (holomorphic) $\P$- vector bundle of rank $r$ over the real (complex) $\P$-manifold $\M = (M, \S)$ . Then it admits local trivializations $\varphi_i: \E|_{U_i} \rightarrow \rund{\S|_{U_i}}^{\oplus r}$ being $\S|_{U_i}$-module 
isomorphisms for a suitable open cover $M = \bigcup_{i \in I} U_i$ , $U_i \subset M$ open, $i \in I$ , together with $\P$- transition functions $A_{i j} \in GL\rund{r, \S\rund{U_i \cap U_j}}$ such that

\[
\begin{array}{ccc}
 & \E|_{U_i \cap U_j} & \\
{}^{\left.\varphi_i\right|_{U_i \cap U_j}} \swarrow & \circlearrowleft & \searrow^{\left.\varphi_j\right|_{U_i \cap U_j}} \phantom{12} \\
\rund{\S|_{U_i \cap U_j}}^{\oplus r} & \mathop{\longrightarrow}\limits_{A_{i j}} & \rund{\S|_{U_i \cap U_j}}^{\oplus r}
\end{array}
\]

for all $i, j \in I$ . $\E^\#$ then is given by local trivializations $U_i \times \rz^r$ (~$U_i \times \cz^r$ ) , $i \in I$ , together with the ordinary transition functions

\[
A_{i j}^\# \in GL\rund{r, \C^\infty\rund{U_i \cap U_j}} \, ( \, A_{i j}^\# \in GL\rund{r, \O\rund{U_i \cap U_j}} \, )  \, .
\]

Again the kernel of ${}^\# : \E \rightarrow \Gamma^{\infty}(\diamondsuit, E)$ ( ${}^\# : \E \rightarrow \Gamma^{\hol}(\diamondsuit, E)$ ) is $\I \E$ (~$\I^\cz \E$ ), and we have a canonical sheaf isomorphism \\
$\Gamma^{\infty}(\diamondsuit , E) \otimes \I^n \simeq \I^n \E$ ( $\Gamma^{\hol}(\diamondsuit, E) \otimes \rund{\I^\cz}^n \simeq \rund{\I^\cz}^n \E$ ) whenever $n \in \nz$ such that $\I^{n + 1} = 0$ .

\item[(iii)] If in addition $\Phi = \rund{\varphi, \rund{\phi_U}_{U \subset M \text{ open } }} : \N \hookrightarrow_{\P} \M$ is a $\P$-imbedding of the real (complex) $\P$-manifold $\N = (N, \T)$ into $\M$ we get canonical maps 
$\E(U) \rightarrow \E|_N\rund{\varphi^{- 1}(U)} \rightarrow \E|_{\N}\rund{\varphi^{- 1}(U)}$ for all $U \subset M$ open, which are called the canonical restrictions, respecting the local trivializations $\rund{\T|_{\varphi^{- 1}\rund{U_i}} }^r$ of $\E|_{\N}$ with the 
$\P$- transition functions

\[
\phi_{U_i \cap U_j} \rund{A_{i j}} \in GL\rund{r, \T\rund{\varphi^{- 1}\rund{U_i} \cap \varphi^{- 1}\rund{U_j}} } \, .
\]

\item[(iv)] Now let the real (holomorphic) $\P$- vector bundles $\E$ and $\F$ over the real (complex) $\P$-manifold $\M = (M, \S)$ be given by local trivializations $U_i \times \rz^r$ ( $U_i \times \cz^r$ ) , $i \in I$ , together with the $\P$- transition functions
$A_{i j} \in GL\rund{r, \S\rund{U_i \cap U_j}}$ resp. $B_{i j} \in GL\rund{s, \S\rund{U_i \cap U_j}}$ then $\E \otimes \F$ is given by the $\P$- transition functions

\[
A_{i j} \otimes B_{i j} \in GL\rund{r s, \S\rund{U_i \cap U_j}} \, .
\]

Again if in addition $\Phi : \N \hookrightarrow_{\P} \M$ is a $\P$-imbedding of the real (complex) $\P$-manifold $\N = (N, \T)$ into $\M$ we see that \\
$(\E \otimes \F)|_{\N} = \E|_{\N} \otimes \F|_{\N}$ .

\item[(v)] One can show that if $\E$ is a real $\P$- vector bundle over the real smooth manifold $M$ then $\E \simeq \E^\#$ as $\P$- vector bundles. This is again rigidity under local deformations.
\end{itemize}

\begin{examples}
\end{examples}

Let $\M = (M, \S)$ be a real (complex) $\P$-manifold of dimension $n$ given by local $\P$-charts $U_i \subset \rz^n$ ( $U_i \subset \cz^n$ ) open, $i \in I$ , forming an open cover of $M$ , and glueing data 
$\Phi_{i j}: U_i \mathop{\supset}\limits_{ \text{ open } } U_{i j} \rightarrow_{\P} U_{j i} \mathop{\subset}\limits_{ \text{ open } } U_j$~. Then the tangent bundle $T \M$ and the cotangent bundle $T^* \M$ of $\M$ are the real (holomorphic) $\P$- vector bundles 
on $\M$ of rank $n$ given by the local trivializations $\rund{\S|_{U_i}}^{\oplus n}$~, $i \in I$ , with transition functions $D \Phi_{i j}$ resp. $\rund{D \Phi_{i j}}^{- 1}$~, where the Jacobian $D \Phi_{i j} \in GL\rund{n, \S\rund{U_{i j}} }$ is taken componentwise 
from the tuple \\
$\rund{f_1, \dots, f_n} \in \rund{\S\rund{U_{i j}} }^{\oplus n}$ associated to $\Phi_{i j}$ , $i, j \in I$ . So clearly \\
$\rund{T \M}^\# = T M$ and $\rund{T^* \M}^\# = T^* M$ .

\begin{lemma} \label{globsect}
Let $\M$ be a complex $\P$-manifold and $\E$ be a holomorphic $\P$- vector bundle over $\M$ with body $E \rightarrow M$ . Let $d := \dim H^0(E) < \infty$ . Then

\[
d \leq \dim H^0(\E) \leq d \dim \P  \, ,
\]

and equivalent are

\begin{itemize}
\item[(i)] $\dim H^0(\E) = d \dim \P$ ,
\item[(ii)] there exist $f_1, \dots, f_d \in H^0(\E)$ such that $\rund{f_1^\#, \dots, f_d^\#}$ is a basis of $H^0(E)$ ,
\item[(iii)] $H^0(\E)$ is a free module over $\P^\cz$ of rank $d$ .
\end{itemize}

Furthermore, if $f_1, \dots, f_d \in H^0(\E)$ such that $\rund{f_1^\#, \dots, f_d^\#}$ is a basis of $H^0(E)$ then $\rund{f_1, \dots, f_d}$ is a $\P^\cz$-basis of $H^0(\E)$ .
\end{lemma}

{\it Proof:} Let $\I \lhd \P$ denote the unique maximal ideal of $\P$ . The first inequality is of course trivial if $\I = 0$ .  For $\I \not= 0$ let $N' \in \nz$ be maximal such that $\I^{N'} \not= 0$ . Then 
$H^0(E) \otimes \rund{\I^\cz}^{N'} = \rund{\I^\cz}^{N'} H^0(\E) \sqsubset H^0(\E)$ , which proves the first inequality. \\

The second inequality, the implication (i) $\Rightarrow$ (ii) and the last statement will be proven by induction on $N \in \nz \setminus \{0\}$ such that $\I^N = 0$ . If $N = 1$ then $\I = 0$ , and all statements are trivial.

Now assume $\I^{N + 1} = 0$ . Then again define $\Q := \P \left/ \I^N \right.$ with unique maximal ideal $\J := \I \left/ \I^N \right.$ having $\J^N = 0$ , and let ${}^\natural : \P \rightarrow \Q$ be the canonical projection. Let 
$\E^\natural := \E \left/ \rund{\I^\cz}^N \E \right.$ , which is a holomorphic $\Q$- vector bundle over $\M^\natural$ of the same rank as $\E$ , and let

\[
{}^\natural : H^0(\E) \rightarrow H^0\rund{\E^\natural}
\]

be the linear map induced by the canonical sheaf projection $\E \rightarrow \E^\natural$ . Its kernel is $\rund{\I^\cz}^N H^0(\E) = H^0(E) \otimes \rund{\I^\cz}^N$ . By induction hypothesis $\dim H^0\rund{\E^\natural} \leq d \dim \Q$ , and so

\[
\dim H^0(\E) \leq d \dim \Q + d \dim \I^N = d \dim \P  \, ,
\]

which proves the second inequality. \\

For proving the implication (i) $\Rightarrow$ (ii) assume $\dim H^0(\E) = d \dim \P$ . Then since $\dim \P = \dim \Q + \dim \I^N$ , $\dim\rund{H^0(E) \otimes \rund{\I^\cz}^N} = d \dim \I^N$ and $\dim H^0\rund{\E^\natural} \leq d \dim \Q$ , we see that necessarily 

\[
{}^\natural : H^0(\E) \rightarrow H^0\rund{\E^\natural}
\]

is surjective and $\dim H^0\rund{\E^\natural} = d \dim \Q$ . So by induction hypothesis and surjectivity there exist $f_1, \dots, f_d \in H^0\rund{\E}$ such that $\rund{f_1^\# , \dots, f_d^\#}$ is a basis of $H^0(E)$ , which proves (ii) . \\

For proving the last statement let $f_1, \dots, f_d \in H^0(\E)$ such that $\rund{f_1^\# , \dots, f_d^\#}$ is a basis of $H^0(E)$ . Then by induction hypothesis $\rund{f_1^\natural, \dots, f_d^\natural}$ is a $\Q^\cz$-basis of $H^0\rund{\E^\natural}$ . For 
proving that $\rund{f_1, \dots, f_d}$ spans $H^0(\E)$ over $\P^\cz$ let $F \in H^0(\E)$ . Then there exist $a_1, \dots, a_d \in \P^\cz$ such that

\[
F^\natural = a_1^\natural f_1^\natural + \dots + a_d^\natural f_d^\natural  \, ,
\]

and so

\[
\Delta := F - a_1^\natural f_1^\natural - \dots - a_d^\natural f_d^\natural \in \rund{\I^\cz}^N H^0(\E) = H^0(E) \otimes \rund{\I^\cz}^N \, .
\]

Since $\rund{f_1^\#, \dots, f_d^\#}$ is a basis of $H^0(E)$ we see that there exist $b_1, \dots, b_d \in~\rund{\I^\cz}^N$ such that

\[
\Delta = f_1^\# \otimes b_1 + \dots + f_d^\# \otimes b_d = b_1 f_1 + \dots + b_d f_d  \, ,
\]

and so

\[
F = \rund{a_1 + b_1} f_1 + \dots + \rund{a_d + b_d} f_d  \, .
\]

For proving linear independence let $a_1, \dots, a_d \in \P^\cz$ such that

\[
a_1 f_1 + \dots + a_d f_d = 0 \, .
\]

Then $a_1^\natural f_1^\natural + \dots + a_d^\natural f_d^\natural = 0$ in $H^0\rund{\E^\natural}$ , and so $a_1^\natural = \dots = a_d^\natural = 0$ . Therefore $a_1, \dots, a_d \in \rund{\I^\cz}^N$ , and this means

\[
0 = a_1 f_1 + \dots + a_d f_d = f_1^\# \otimes a_1 + \dots + f_d^\# \otimes a_d \, .
\]

Since $f_1^\#, \dots, f_d^\#$ are linearly independent we get $a_1 = \dots = a_d = 0$ . \\

Now (ii) $\Rightarrow$ (iii) follows from the last statement, and (iii) $\Rightarrow$ (i) is of course trivial. $\Box$

\section{$\mathcal{P}$-lattices and automorphic forms} \label{lat and aut}

Let $G$ be a real {\sc Lie} group. Then it is in particular a smooth real manifold, and the multiplication on $G$ can be written as a smooth map $m: G \times G \rightarrow G$~. Therefore the multiplication turns the set $G^{\P}$ of all $\P$-points of $G$ into a 
group via $g h := m(g, h)$ for all $g, h \in_{\P} G$ , and clearly ${}^\#: G^{\P} \rightarrow G \, , \, g \mapsto g^\#$ is a group epimorphism. Of course the $\P$-points of $GL(n, \rz)$ are in 1-1 correspondence with $n \times n$ -matrices having entries in 
$\P$ and body (taken componentwise) in $GL(n, \rz)$ , and the product of two $\P$-points of $GL(n, \rz)$ can be computed via ordinary matrix multiplication.

\begin{defin} [ $\P$-lattices]
Let $\Upsilon$ be a subgroup of $G^{\P}$ . $\Upsilon$ is called a $\P$-lattice of $G$ iff
\begin{itemize}
\item[\{i\}] $\Upsilon^\# := \schweif{\left.\gamma^\# \, \right| \, \gamma \in \Upsilon} \sqsubset G$ is an ordinary lattice, called the body of $\Upsilon$ , and
\item[\{ii\}] ${}^\# : \Upsilon \rightarrow \Upsilon^\# \, , \, \gamma \mapsto \gamma^\#$ is bijective and so automatically an isomorphism.
\end{itemize}
\end{defin}

Obviously a $\P$-lattice $\Upsilon$ is nothing but a local deformation over the algebra $\P$ of the natural embedding $\Gamma := \Upsilon^\# \hookrightarrow G$ as a group homomorphism. Of course given a $\P$-lattice $\Upsilon$ of $G$ with body 
$\Gamma \sqsubset G$ and $g \in_{\P} G$ with $g^\# = 1$ we get another $\P$-lattice $g \Upsilon g^{- 1}$ of $G$ with same body $\Gamma$ . The set of all $\P$-lattices of $G$ of the form $g \Upsilon g^{- 1}$ of $G$ , $g \in_{\P} G$ , $g^\# = 1$ , is called the 
conjugacy class of $\Upsilon$ . In the case $\I^2 = 0$~, where $\I \lhd \P$ denotes the unique maximal ideal of $\P$ , given an ordinary lattice $\Gamma$ the conjugacy classes of $\P$-lattices $\Upsilon$ with body $\Gamma$ 
are in 1-1 correspondence with $H^1\rund{\Gamma, \g} \otimes \I$ , $\Gamma$ acting on the {\sc Lie} algebra $\g$ of $G$ by $\Ad$ , see for example \cite{Ragh}~.

\begin{lemma} \label{powers}
\item[(i)] Let $\Upsilon$ be a $\P$-lattice of $G$ , $\gamma \in_{\P} \Upsilon$ and $n \in \nz \setminus \{0\}$ such that $\rund{\gamma^\#}^n = 1$~. Then $\gamma^n = 1$ .
\item[(ii)] Let $g \in_{\P} G$ with $g^n = 1$ and $g^\# \in Z(G)$ . Then $g = g^\#$ .
\end{lemma}

{\it Proof:} (i) Obviously $\gamma^n \in \Upsilon$ having $\rund{\gamma^n}^\# = 1$ . Therefore $\gamma^n = 1$ by property \{ii\} . $\Box$ \\

(ii) Clearly $\rund{g^\#}^n = \rund{g^n}^\# = 1$ . And so it suffices to show that \\
$\varphi: G \rightarrow G \, , \, h \mapsto h^n$ is a local diffeomorphism at $g^\#$ . But, since $g^\# \in Z(G)$ one can easily compute

\[
D \varphi \rund{g^\#} = n \, \rund{D t_{g^\#}(1)}^{- 1}  \, ,
\]

which is invertible since the translation $t_{g^\#} : G \rightarrow G$ with $g^\#$ is a diffeomorphism. $\Box$ \\

Let $H := \{z \in \cz \, | \, \Im z > 0\}$ be the usual upper half plane, and from now on let $G := SL(2, \rz)$ . Then $G$ acts on $H$ via {\sc Möbius} transformations

\[
g z := \frac{a z + b}{c z + d} \, , \, g = \rund{\begin{array}{cc} a & b \\
c & d \end{array}} \, ,
\]

more precisely we have a group epimorphism $\overline{\phantom{1}} : G \rightarrow \Aut(H)$ with kernel $\{\pm 1\} = Z(G)$ . The action of $G$ on $H$ induces a morphism of ringed spaces

\[
\rund{G \times H, \pr_1^* \rund{\C^\infty_G}^\cz \hat\otimes \pr_2^* \O_H} \rightarrow \rund{H, \O_H} \, ,
\]

and therefore a group homomorphism $G^{\P} \rightarrow \{\P\text{-automorphisms of } H\}$ respecting $\phantom{1}^\#$ with kernel $\{\pm 1\}$ , which is no longer surjective if $\P \not= \rz$ . If $g \in_{\P} G$ then $g$ as a $\P$-automorphism of $H$ is 
given by $g z \in \O(H) \otimes \P^\cz$ . For all $U \subset H$ open and $f \in \O(U) \otimes \P^\cz$ denote by $f(g z) \in \O\rund{g^{- 1} U} \otimes \P^\cz$ the pullback of $f$ under $g$ . If $U$ is invariant under $g^\#$ then we say $f$ is $g$-invariant if 
and only if $f(gz) = f$ . \\

Let $k \in \nz$ be fixed for the rest of the article, and let $j \in \C^\infty(G)^\cz \hat\otimes \O(B)$ be given by

\[
j(g, z) := \frac{1}{c z + d} \, , \, g = \rund{\begin{array}{cc} a & b \\
c & d \end{array}} \, .
\]

Then $j$ fulfills the cocycle property $j(g h, z) = j(g, h z) j(h, z)$ , and an easy computation shows that $j(g, z)^2 = g'(z)$ , $g$ regarded as an automorphism of $H$ . \\

This gives a right-action of $G$ on $\O(H)$ by

\[
| : \O(H) \rightarrow \C^\infty(G)^\cz \hat\otimes \O(H) \, , \, f|_g(z) := f(g z) j(g, z)^k \, ,
\]

and this induces an action of $G^{\P}$ on $\O(H) \otimes \P^\cz$ given by

\[
|_g : \O(H) \otimes \P^\cz \rightarrow \O(H) \otimes \P^\cz \, , \, f|_g(z) := f(g z) j(g, z)^k
\]

for all $g \in_{\P} G$ , or more precisely for all $U \subset H$ this gives a $\P^\cz$-linear map

\[
|_g : \O(U) \otimes \P^\cz \rightarrow \O\rund{g^{- 1} U} \otimes \P^\cz \, , \, f|_g(z) := f(g z) j(g, z)^k  \, .
\]

From now on let $\Upsilon$ be a fixed $\P$-lattice in $G$ with body $\Gamma := \Upsilon^\#$ .

\pagebreak

\begin{examples} \label{paramlatt}
\end{examples}

\begin{itemize}

\item[(i)] Let $\Gamma := SL(2, \zz)$ . Then $\Gamma$ is the free group generated by the two matrices $R := \rund{\begin{array}{cc}
0 & 1 \\
- 1 & - 1
\end{array}}$ and $S := \rund{\begin{array}{cc}
0 & 1 \\
- 1 & 0
\end{array}}$ modulo the relations $R^3 = S^4 = 1$ . One can easily compute that the equations $g^3 = h^4 = 1$ define a smooth submanifold $M$ of $G^2$ of dimension $4$ near the point $(R, S)$ and that the map

\[
\varphi: G \rightarrow M \, , \, g \mapsto \rund{g^{- 1} R g, g^{- 1} S g}
\]

has injective differential at $g = 1$ . So take any smooth submanifold $M' \subset M$ of dimension $1$ such that $T_{(R, S)} M = T_{(R, S)} M' \oplus \Im D \varphi(1)$~. Then obviously the conjugacy classes of $\P$-lattices $\Upsilon$ with 
$\Upsilon^\# = \Gamma$ are in 1-1 correspondence with $\P$-points $x \in_{\P} M'$ having $x^\# = (R, S)$ , and so via a local chart of $M'$ at $(R, S)$ with $\I$ . In particular

\[
H^1(\Gamma, \g) = \left.T_{(R, S)} M \right/ \Im D \varphi(1) \simeq \g \left/ \rund{\z_\g(R) + \z_\g(S)} \right.
\]

has dimension $1$ .

\item[(ii)] Let $X$ be a compact {\sc Riemann} surface of genus $g$ , $s_1, \dots, s_m \in X$ , \\
$3 g + m \geq 3$ . Then the universal covering of $X \setminus \{s_1, \dots, s_m\}$ is isomorphic to $H$ , and by \cite{Nat} one can write 
$X \setminus \{s_1, \dots, s_m\} =~\Gamma~\backslash~H$~, where $\Gamma \subset G$ is a lattice without elliptic elements having $- 1 \notin~\Gamma$~, it is the free group generated by some hyperbolic elements 
$A_1, B_1, \dots, A_g, B_g \in~G$ and parabolic elements $C_1, \dots, C_n \in G$ modulo the single relation

\[
\eckig{A_1, B_1} \dots \eckig{A_g, B_g} C_1 \dots C_m = 1 \, .
\]

Then by exactly the same method as in (i) one obtains a 1-1 correspondence between the conjugacy classes of $\P$-lattices $\Upsilon$ with $\Upsilon^\# = \Gamma$ and $\P$-points $x$ of a suitable $\rund{6 (g - 1) + 3 m}$-dimensional smooth 
submanifold of $G^{2 g + m}$ having $x^\# = \rund{A_1, B_1, \dots, A_g, B_g, C_1, \dots, C_m}$ and so with $\I^{\oplus \rund{6 (g - 1) + 3 m}}$ . In particular

\[
\dim H^1(\Gamma, \g) = 6 (g - 1) + 3 m \, .
\]

\end{itemize}

For defining automorphic forms with respect to the $\P$-lattice $\Upsilon$ we need some more informtion about the behaviour of $\Gamma$ being an ordinary 
lattice in $G$~. Let $\pz^1$ denote the {\sc Riemann} sphere, on which $SL(2, \cz)$ acts via {\sc Möbius} transformations. For $z \in H \cup \partial_{\pz^1} H$ denote by $\overline z := \Gamma z$ the image of $z$ under the canonical projection 
$H \rightarrow \Gamma \backslash H$ resp. $\partial_{\pz^1} H \rightarrow \Gamma \backslash \partial_{\pz^1} H$ .

\begin{defin}
\item[(i)] An element $\overline{z_0} \in \Gamma \backslash H$ , $z_0 \in H$ , is called regular iff $\overline \Gamma^{z_0} = \{\id\}$~,
\item[(ii)] an element $\overline{z_0} \in \Gamma \backslash H$ , $z_0 \in H$ , is called elliptic iff $\overline \Gamma^{z_0} \not= \{\id\}$~, and Finally,
\item[(iii)] an element $\overline{z_0} \in \Gamma \backslash \partial_{\pz^1} H$ , $z_0 \in \partial_{\pz^1} H$ , is called a cusp of $\Gamma \backslash H$ iff $\overline \Gamma \cap \overline{P^{z_0}} \not= \{\id\}$ , where $P^{z_0} \sqsubset G$ denotes the 
parabolic subgroup associated to $z_0$ .
\end{defin}

It is a well known fact that there exist always only finitely many elliptic points in $\Gamma \backslash H$ ,
$\Gamma \backslash H$ has always only finitely many cusps, and the quotient $\Gamma \backslash H$ can be compactified as a topological space by adding the cusps of $\Gamma \backslash H$ . This can for example be deduced from theorem 0.6 in 
\cite{GarlRagh} . \\

Since $G$ acts transitively on $H$ we see that for each $z_0 \in H$ there exists $g \in G$ such that $g i = z_0$ , and therefore $G^{z_0} = g K g^{- 1}$ , where $K := G^i \simeq \rz / \zz$ is a maximal compact subgroup of $G$ . Therefore if $z_0$ is an 
elliptic point of $\Gamma \backslash H$ then $\Gamma^{z_0} \sqsubset G^{z_0}$ and $\overline{\Gamma}^{z_0} \sqsubset \overline G^{z_0}$ are finite non-trivial cyclic groups. $\ord \ \overline{\Gamma}^{z_0} \in \nz$ is called the period of 
$\overline{z_0}$ . \\

Since Furthermore, $G$ acts transitively on the boundary $\partial_{\pz^1} H$ of $H$ we see that for each $z_0 \in \partial_{\pz^1} H$ there exists an element $g \in G$ such that $g(\infty) = z_0$ , and so $P^{z_0} = g P^\infty g^{- 1}$ . Recall that 
$P^\infty \simeq \rz$ is the one-parameter-subgroup generated by $\chi_0 \in \g$ , $\g$ being the {\sc Lie} algebra of $G$ , with

\[
\chi_0 = \rund{\begin{array}{cc} 0 & 1 \\
0 & 0
\end{array}} \, .
\]

Therefore if $z_0$ is a cusp of $\Gamma \backslash H$ then $\overline \Gamma \cap \overline{P^{z_0}}$ is infinite cyclic, and one can always choose $g \in G$ such that in addition 
$g^{- 1} \overline \Gamma g \cap \overline{P^\infty} = \spitz{\overline{g_0}}$ , where

\[
g_0 := \rund{\begin{array}{cc} 1 & 1 \\
0 & 1
\end{array}} = \exp\rund{\chi_0} \, .
\]

\begin{lemma} \label{paramparab}
\item[(i)] Let $g \in_{\P} G$ such that $g^\# = g_0$ . Then there exists a unique $\chi \in \g \otimes \P$ such that $\chi^\# = \chi_0$ and $g = \exp(\chi)$ .
\item[(ii)] Let $\chi \in \g \otimes \P$ with body $\chi^\# = \chi_0$ . Then there exists a $\P$-automorphism $\Omega: H \rightarrow_{\P} H$ such that $\Omega^\# = \id$ and

\[
\begin{array}{ccc}
\phantom{123456789012345678,} H & \mathop{\longrightarrow}\limits^{\Omega} & H \phantom{1234567} \\
\exp\rund{t \chi_0} : z \mapsto z + t \downarrow & \circlearrowleft & \downarrow \exp(t \chi) \\
\phantom{123456789012345678,} H & \mathop{\longrightarrow}\limits_{\Omega} & H \phantom{1234567}
\end{array}
\]

for all $t \in \rz$ . All other $\P$-automorphisms with this property are given by $z \mapsto \Omega\rund{z + a}$ where $a \in \I^\cz$ .
\end{lemma}

{\it Proof:} (i) For proving this statement it suffices to show that $\exp$ is a local diffeomorphism at $\chi_0$ . We use theorem 1.7 of chapter II section 1.4 in \cite{Helga}~, which says the following:

\begin{quote}
Let $G$ be a {\sc Lie} group with {\sc Lie} algebra $\g$ . The exponential mapping of the manifold $\g$ into $G$ has the differential

\[
D \exp_X = D \rund{l_{\exp X}}_e \circ \frac{1 - e^{- \ad_X}}{\ad_X}  \phantom{12345} (X \in \g) \, .
\]

As usual, $\g$ is here identified with the tangent space $\g_X$ .
\end{quote}

Hereby $e$ denotes the unit element of the {\sc Lie} group $G$ , $l_g$ denotes the left translation on $G$ with an element $g \in G$ , 	and $\exp$ is used as a local chart of $G$ at $e$ . \\

Since $\chi_0$ is nilpotent in $\g$ we see that also $\ad_{\chi_0} \in \End(\g)$ is nilpotent, and so is $\frac{1 - e^{- \ad_{\chi_0}} }{\ad_{\chi_0}} - 1 \in \End(\g)$ . Since $l_{\exp{\chi_0}} : G \rightarrow G$ is a diffeomorphism, we obtain the desired 
result applying the theorem with $X := \chi_0$ . $\Box$ \\

(ii) Since $\chi^\# = \chi_0$ an easy calculation shows that $\chi$ is nilpotent as a matrix with entries in $\P$ . Therefore since in addition $\exp\rund{t \chi_0}$ is upper triangle we see that $\omega:= \exp(t \chi) \, i \in \P^\cz [t]$ with body 
$\omega^\# = \exp\rund{t \chi^\#} \, i = t + i$ . Now let $\Omega := \omega(t - i) \in \P^\cz [t]$ . Then $\Omega^\# = t$ , and so $\Omega$ can be regarded as a $\P$-automorphism $\Omega: H \rightarrow H$ having $\Omega^\# = \id$ .

Since everything in the diagramme is given by tuples of holomorphic functions on $H$ it suffices to prove its commutativity on the non discrete subset \\
$\rz + i \subset H$ . So let $t, u \in \rz$ . Then

\begin{eqnarray*}
\rund{\Omega \circ \exp\rund{t \chi_0}} (u + i) &=& \Omega(u + i + t) \\
&=& \omega(u + t) \\
&=& \exp(t \chi) \exp(u \chi) i \\
&=& \exp(t \chi) \omega(u) \\
&=& \rund{\exp(t \chi) \circ \Omega}(u + i)  \, .
\end{eqnarray*}

Now let $\widetilde \Omega : H \rightarrow_{\P} H$ be another $\P$-automorphism. Then $\widetilde \Omega$ has the same properties iff $\Omega^{- 1} \circ \widetilde \Omega$ is a $\P$-automorphism with body 
$\id$ and commuting with all translations $H \rightarrow H \, , \, z \mapsto z + t$ , $t \in \rz$ , iff $\rund{\Omega^{- 1} \circ \widetilde \Omega} (z) = z + a$ with some $a \in \O(H) \otimes \I^\cz$ and invariant under the translations 
$H \rightarrow H \, , \, z \mapsto z + t$ , $t \in \rz$ , and therefore constant. $\Box$ \\

\begin{defin}

Let $z_0 \in \partial_{\pz^1} H$ , $\gamma \in_{\P} G$ such that $\gamma^\# \in P^{z_0} \setminus \{\pm 1\}$ , $U \subset H$ open and $\gamma^\#$-invariant and Finally, $f \in \O(U) \otimes \P^\cz$ such that $f|_\gamma = f$ or $f$ $\gamma$-invariant (which 
is nothing but $f|_\gamma = f$ for $k = 0$ ). Let $g \in G$ such that $g \infty = z_0$ and $g_0 = g^{- 1} \gamma^\# g$ . Furthermore, assume that there exists $c > 0$ such that

\[
\{\Im z > c\} \subset g^{- 1} U \, ,
\]

which is $g_0$-invariant. Let $\chi$ and $\Omega$ be given by lemma \ref{paramparab} taken \\
$\widetilde{g_0} := g^{- 1} \gamma^\# g$ instead of $g$ having body $g_0$ . If we define $\left.f|_g \right|_\Omega$ as

\[
\left.f|_g \right|_\Omega (z) := f|_g \rund{\Omega z} \Omega'(z)^{k / 2}
\]

then we see that

\[
\left.f|_g \right|_\Omega (z) = \left.\left.f|_g \right|_\Omega \right|_{g_0} (z) = \left.f|_g \right|_\Omega (z + 1) \, .
\]

Now $f$ is called bounded (vanishing) at $z_0$ iff $\left.f|_g \right|_\Omega (z)$ is bounded, and therefore converging, (resp. vanishing) for $\Im z \leadsto \infty$ .
\end{defin}

Observe that $\rund{\Omega'}^r$ is well defined for all $r \in \rz$ since $\Omega^\# = \id$ , and so $\rund{\Omega'}^\# = 1$ . Clearly the definition does not depend on the choice of $g$ and $\Omega$ because $g$ is uniquely determined up to $\pm 1$ , and
if $\widetilde \Omega$ is another choice for $\Omega$ then $\widetilde \Omega(z) = \Omega(z + a)$ by lemma \ref{paramparab} with some $a \in \I^\cz$ . Therefore $\left.f|_g \right|_{\widetilde \Omega} (z) = \left.f|_g \right|_\Omega (z + a)$ . \\

\begin{defin}[automorphic and cusp forms for $\Upsilon$ ]
Let \\
$f \in \O(H) \otimes \P^\cz$ . $f$ is called an automorphic (cusp) form for $\Upsilon$ of weight $k$ iff

\begin{itemize}
\item[(i)] $f|_\gamma = f$ for all $\gamma \in_{\P} \Upsilon$ ,
\item[(ii)] $f$ is bounded (vanishing) at all cusps $\overline{z_0} \in \Gamma \left\backslash \partial_{\pz^1} H \right.$ of $\Gamma \backslash H$ .
\end{itemize}

The space of automorphic (cusp) forms for $\Upsilon$ of weight $k$ is denoted by $M_k(\Upsilon)$ (resp. $S_k(\Upsilon)$ ) . We have $S_k(\Upsilon) \sqsubset M_k(\Upsilon) \sqsubset \O(H) \otimes \P^\cz$ as $\P^\cz$-submodules.
\end{defin}

Since $\rund{f|_g}^\# = \left.f^\#\right|_{g^\#}$ for all $f \in \O(H) \otimes \P^\cz$ and $g \in_{\P} G$ we see that $M_k(\Upsilon)^\# = M_k(\Gamma)$ and $S_k(\Upsilon)^\# = S_k(\Gamma)$ . Using lemma \ref{powers} we observe that 
$- 1 \in \Gamma \Leftrightarrow - 1 \in \Upsilon$ , and so in this case $M_k(\Upsilon) = 0$ if $2 \not| k$ .

\begin{theorem}

\begin{eqnarray*}
&& \dim M_k(\Gamma) \leq \dim M_k(\Upsilon) \leq \dim M_k(\Upsilon) \dim \P \\
&& ( \, \dim S_k(\Gamma) \leq \dim S_k(\Upsilon) \leq \dim S_k(\Upsilon) \dim \P \, ) \, ,
\end{eqnarray*}

and equivalent are

\begin{itemize}
\item[(i)] $\dim M_k(\Upsilon) = \dim M_k(\Gamma) \dim \P$ ( $\dim S_k(\Upsilon) = \dim S_k(\Gamma) \dim \P$ ),
\item[(ii)] there exist $f_1, \dots, f_r \in M_k(\Upsilon)$ ( $f_1, \dots, f_r \in S_k(\Upsilon)$ ) such that $\rund{f_1^\#, \dots, f_r^\#}$ is a basis of $M_k(\Gamma)$ ( $S_k(\Gamma)$ ),
\item[(iii)] $M_k(\Upsilon)$ ( $S_k(\Upsilon)$ ) is a free module over $\P^\cz$ of rank $\dim M_k(\Gamma)$ (~$\dim S_k(\Gamma)$ ).
\end{itemize}

Furthermore, if $f_1, \dots, f_r \in M_k(\Upsilon)$ ( $f_1, \dots, f_r \in S_k(\Upsilon)$ ) such that $\rund{f_1^\#, \dots, f_r^\#}$ is a basis of $M_k(\Gamma)$ ( $S_k(\Gamma)$ ) then $\rund{f_1, \dots, f_r}$ is a $\P^\cz$-basis of $M_k(\Upsilon)$ ( $S_k(\Upsilon)$ ).
\end{theorem}

{\it Proof:} This is a corollary of lemma \ref{globsect} since in section \ref{quot} and \ref{sect main} we will show that $M_k(\Gamma)$ and $S_k(\Gamma)$ are the spaces of global sections for certain holomorphic $\P$- line bundles on 
$\Upsilon \backslash H \cup \{ \text{ cusps of } \Gamma \backslash H\}$ as $\P$- {\sc Riemann} surface. $\Box$ \\

The aim of this article is now to prove that in almost all cases we have an isomorphism $M_k(\Upsilon) \simeq M_k(\Gamma) \otimes \P^\cz$ mapping $S_k(\Upsilon)$ to $S_k(\Gamma) \otimes \P^\cz$ .

\section{The quotient as a $\mathcal{P}$- {\sc Riemann} surface} \label{quot}

It is a well known fact that there exists a structure of a compact {\sc Riemann} surface on $X := \Gamma \backslash H \cup \schweif{ \text{ cusps of } \Gamma \backslash H}$ such that the subsheaf of $\O_H$ of $\Gamma$-invariant functions is the pullback 
of $\left.\O_X \right|_{\Gamma \backslash H}$ under the canonical projection $\pi: H \rightarrow \Gamma \backslash H \hookrightarrow X$ . Now we will construct a $\P$- {\sc Riemann} surface \\
$\X = (X, \S)$ such that the subsheaf of $\O_H \otimes \P^\cz$ of $\Upsilon$-invariant functions is precisely the pullback of $\S |_{\Gamma \backslash H}$ under the canonical projection \\
$\pi: H \rightarrow \Gamma \backslash H$ . For this purpose define the sheaf $\S$ of $\P^\cz$-algebras on $X$ as

\begin{eqnarray*}
&& \S(V) := \left\{f \in \O\rund{\pi^{- 1}(V)} \otimes \P^\cz \phantom{12} \Upsilon\text{-invariant and bounded at } \right. \\
&& \left. \phantom{1234 \pi^{- 1} \P^\cz} \text{all cusps } \overline{z_0} \in V \text{ of } \Gamma \backslash H\right\}
\end{eqnarray*}

for all $V \subset X$ open. Recall that a function $f \in \O\rund{\pi^{- 1}(V)} \otimes \P^\cz$ is called $\Upsilon$-invariant iff $f(\gamma z) = f$ for all $\gamma \in \Upsilon$ . Now one has to show that locally $\S \simeq \O_X \otimes \P^\cz$ . We will do this 
giving local $\P$-charts for $\X$ . On may define the $\P$- {\sc Riemann} surface $\Upsilon \backslash H := \rund{\Gamma \backslash H, \S|_{\Gamma \backslash H}}$ to be the quotient of $H$ by $\Upsilon$ and the $\P$-morphism 
$\Pi := \rund{\pi, \rund{i_V}_{V \subset \Gamma \backslash H \text{ open }} }$ from $H$ to $\Upsilon \backslash H$ as the canonical $\P$-projection, where $i_V : \S(V) \hookrightarrow \O\rund{\pi^{- 1}(V)} \otimes \P^\cz$ , $V \subset \Gamma \backslash H$ 
open, denote the canonical inclusions. \\

Let $\F$ be the $\S$- sheaf module on $X$ given by

\begin{eqnarray*}
&& \F(V) := \left\{f \in \O\rund{\pi^{- 1}(V)} \otimes \P^\cz \, \right| \, f|_\gamma = f \text{ for all } \gamma \in \Upsilon \text{ and } \\
&& \left. \phantom{1234 \pi^{- 1} \P^\cz} f \text{ bounded at all cusps } \overline{z_0} \in V \text{ of } \Gamma \backslash H\right\}
\end{eqnarray*}

for all $V \subset X$ open. If $- 1 \in \Gamma$ and $k$ is odd then of course $\F = 0$ , and in the cases where either $k$ is even or $k$ is odd and $- 1 \notin \Gamma$ we will show that $\F$ is a holomorphic $\P$- line bundle over $\X$ . For this 
purpose we have to give local trivializations $\F \simeq \S$ as $\S$- sheaf modules. Then obviously $M_k(\Upsilon) = H^0(\F)$~. \\

For $g \in G$ let $\pi_g : H \rightarrow \spitz{g} \backslash H$ denote the canonical projection. \\

{\it At regular points of $\Gamma \backslash H$ }:

\begin{quote}
Let $\overline{z_0} \in \Gamma \backslash H$ , $z_0 \in H$ , be regular. Then there exists an open neighbourhood $U \subset H$ of $z_0$ auch that $\gamma U \cap U = \emptyset$ for all \\
$\gamma \in \overline\Gamma \setminus \{1\}$ , and so

\[
\pi |_U : U \rightarrow \pi(U) \mathop{\subset}\limits_{ \text{ open } } X
\]

is biholomorphic. Its inverse is a local chart of $X$ at $\overline{z_0}$ . \\

For giving a local $\P$-chart of $\X$ at $\overline{z_0}$ we will show that locally $\Pi$ is an isomorphism at $z_0$ . Indeed the restriction of $\Pi$ to $U$ is given by $\Pi|_U = \rund{\pi|_U, \rund{|_{\pi^{- 1}(V) \cap U} }_{V \subset \pi(U) \text{ open } }}$ from 
$\rund{U, \O_U \otimes \P^\cz}$ to $\rund{\pi(U), \S|_{\pi(U)} }$ as ringed spaces, where for all $V \in \pi(U)$ open

\begin{eqnarray*}
&& |_{\pi^{- 1}(V) \cap U} : \S(V) = \schweif{f \in \O\rund{\pi^{- 1}(V)} \otimes \P^\cz \phantom{12} \Upsilon\text{-invariant } } \\
&& \phantom{1234567890123} \rightarrow \O\rund{\pi^{- 1}(V) \cap U} \otimes \P^\cz
\end{eqnarray*}

simply denotes the restriction map. It is indeed an isomorphism of $\P^\cz$-algebras since $\pi^{- 1}(V) = \dot \bigcup_{\gamma \in \overline \Gamma} \, \gamma \rund{\pi^{- 1}(V) \cap U}$ for all $V \subset \pi(U)$ open. So $\rund{\Pi|_U}^{- 1}$ is a local 
$\P$-chart for $\X$ . \\

For giving a local trivialization of $\F$ at $z_0$ identify the ringed spaces $\rund{U, \O_U \otimes \P^\cz}$ and $\rund{\pi(U), \S|_{\pi(U)} }$ via $\Pi|_U$ . Then we see by the same argument that the restriction maps

\[
|_{\pi^{- 1}(V) \cap U} : \F(V) \rightarrow \O\rund{\pi^{- 1}(V) \cap U} \otimes \P^\cz \, ,
\]

$V \subset \pi(U)$ open, give an isomorphism of the $\S|_{\pi(U)}$- sheaf modules $\left.\F \right|_{\pi(U)}$ and $\O_U \otimes \P^\cz$ .
\end{quote}

{\it At elliptic points of $\Gamma \backslash H$ }:

\begin{quote}
Let $\overline{z_0} \in \Gamma \backslash H$ , $z_0 \in H$ , be elliptic of period $n$ and $g \in G$ such that $g i = z_0$ . Then $\Gamma^{z_0} = \spitz{\gamma_0} \subset G^{z_0} = g K g^{- 1}$ for some $\gamma_0 \in \Gamma^{z_0}$ .

Let $\overline \pi : \spitz{\gamma_0} \backslash H \rightarrow \Gamma \backslash H$ denote the canonical projection. Then $\pi = \overline \pi \circ \pi_{\gamma_0}$ . Now we choose

\[
c: B := \schweif{\left.z \in \cz \frac{}{} \, \right| \, \abs{z} < 1} \rightarrow H \, , \, z \mapsto - i \, \frac{z + 1}{z - 1}
\]

as a {\sc Cayley} transform with $c(0) = i$ . It is clearly biholomorphic, and $c^{- 1} g^{- 1} \gamma_0 g c$ fixes $0$ as an automorphism of $B$ and so is given by multiplication with a suitable $\eta \in U(1)$ of order $n$ . Therefore

\[
\varphi: B \mathop{\longrightarrow}^{\sqrt[n]{\phantom{1}} } \left. \spitz{\eta} \right\backslash B \mathop{\longrightarrow}^c \left. \spitz{g^{- 1} \gamma_0 g} \right\backslash H 
\mathop{\longrightarrow}^g \left. \spitz{\gamma_0} \right\backslash H \mathop{\longrightarrow}^{\overline \pi} \Gamma \backslash H
\]

gives a locally biholomorphic map at $0 \mapsto 0 \mapsto i \mapsto z_0 \mapsto \overline{z_0}$ , and so $\varphi^{- 1}$ is a local chart of $\Gamma \backslash H$ at $\overline{z_0}$ . \\

For giving a local $\P$-chart of $\X$ at $z_0$ define the sheaves $\widetilde \O_B$ , $\widetilde \O_H$ of unital complex algebras and the sheaf $\widetilde \S$ of unital $\P^\cz$-algebras on $\spitz{\eta} \backslash B$ , 
$\left.\spitz{g^{- 1} \gamma_0 g} \right\backslash H$ resp. $\spitz{\gamma_0} \backslash H$ by

\[
\widetilde \O_B(V) := \schweif{\left.f \in \O\rund{\pi_\eta^{- 1}(V)} \, \right| \, f(\eta w) = f} \, ,
\]

\[
\widetilde \O_H(V) := \schweif{\left.f \in \O\rund{\pi_{g^{- 1} \gamma_0 g}^{- 1}(U)} \, \right| \, f\rund{g^{- 1} \gamma_0 g  z} = f}
\]

and

\[
\widetilde \S(V) := \schweif{\left.f \in \O\rund{\pi_{\gamma_0}^{- 1}(V)} \otimes \P^\cz \, \right| \, f\rund{\widetilde{\gamma_0} z} = f}
\]

for all $V \subset \spitz{\eta} \backslash B$ , $V \subset \left.\spitz{g^{- 1} \gamma_0 g} \right\backslash H$ resp. $V \subset \left.\spitz{\gamma_0} \right\backslash H$ , where $\pi_\eta: B \rightarrow \spitz{\eta} \backslash B$ denotes the canonical 
projection and $\widetilde{\gamma_0} \in \Upsilon$ is the unique element such that $\widetilde{\gamma_0}^\# = \gamma_0$ . We will see that $\varphi$ extends to a local isomorphism

\begin{eqnarray*}
&& \Phi: \rund{B, \O_B \otimes \P^\cz} \rightarrow \rund{\spitz{\eta} \backslash B \, , \, \widetilde \O_B \otimes \P^\cz} \\
&& \phantom{12} \rightarrow \rund{\left. \spitz{g^{- 1} \gamma_0 g} \right\backslash H \, , \, \widetilde \O_H \otimes \P^\cz} \rightarrow \rund{\left. \spitz{\gamma_0} \right\backslash H \, , \, \widetilde \S} \rightarrow \X
\end{eqnarray*}

of ringed spaces at $0 \mapsto \overline{z_0}$ , and therefore $\Phi^{- 1}$ is a local $\P$-chart of $\X$ . The first and the second isomorphism are just induced by $\sqrt[n]{\phantom{1}}$ and $c$ .

\begin{lemma}
\item[(i)] There exists a unique $\widetilde{z_0} \in_{\P} H$ such that $\widetilde{\gamma_0} \widetilde{z_0} = \widetilde{z_0}$ , and $\widetilde{z_0}^\# = z_0$ .
\item[(ii)] There exists $\widetilde g \in_{\P} H$ such that ${\widetilde g}^\# = g$ and $\widetilde g \, i = \widetilde z_0$ .
\end{lemma}

{\it Proof:} (i) : Let $E := \{g \in G \text{ elliptic }\} \subset G$ and

\[
M := \{(g, z) \in E \times H \, | \, g z = z\} \, .
\]

Then $M$ is the preimage of $0$ under the smooth map $E \times H \rightarrow H \, , \, (g, z) \mapsto g z - z$ with surjective differential everywhere. One can easily show that $M$ is at the same time the graph of a smooth map $\varphi: E \rightarrow H$ , 
this means it is the preimage of $0$ under the smooth map $E \times H \rightarrow \cz \, , \, (g, z) \mapsto \varphi(g) - z$ with surjective differential everywhere. So 
$\widetilde{\gamma_0} \widetilde{z_0} = \widetilde{z_0} \Leftrightarrow \widetilde z_0 = \varphi\rund{\widetilde{\gamma_0}}$~.~$\Box$ \\

(ii) : One also can easily compute a smooth map $\varphi: H \rightarrow G$ such that $z = \varphi(z) i$ for all $z \in H$ . So take $\widetilde g := \varphi\rund{\widetilde{z_0}} \varphi\rund{z_0}^{- 1} g \in_{\P} G$~. $\Box$ \\

So ${\widetilde g}^{- 1} \widetilde{\gamma_0} \widetilde g \in_{\P} K$ since it fixes $i$ , and $\rund{{\widetilde g}^{- 1} \widetilde{\gamma_0} \widetilde g}^{2 n} = 1$ . Therefore by lemma 
\ref{powers} since $K$ is commutative we have 

\[
{\widetilde g}^{- 1} \widetilde{\gamma_0} \widetilde g = \rund{{\widetilde g}^{- 1} \widetilde{\gamma_0} \widetilde g}^\# = g^{- 1} \gamma_0 g \in K \, .
\]

This gives the commuting diagram

\[
\begin{array}{ccc}
\phantom{123456,} H & \mathop{\longrightarrow}\limits^{\widetilde g} & H \phantom{12.} \\
g^{- 1} \gamma_0 g \downarrow & \circlearrowleft & \downarrow \widetilde{\gamma_0} \\
\phantom{123456,} H & \mathop{\longrightarrow}\limits_{\widetilde g} & H \phantom{12.}
\end{array}
\]

of $\P$-automorphisms inducing the third isomorphism. For the last isomorphism let $U \subset \spitz{\gamma_0} \backslash H$ be an open neighbourhood of $\pi_{\gamma_0}\rund{z_0}$ such that

\[
\left.\overline \pi \right|_{U} : U \rightarrow \overline \pi(U) \mathop{\subset}\limits_{ \text{ open } } X
\]

is biholomorphic. Then for all $V \subset \overline \pi(U)$ open $\pi_{\gamma_0}^{- 1} \rund{\overline \pi^{\, - 1}(V) \cap U} \subset H$ is already $\gamma_0$-invariant, and 

\[
\pi^{- 1}(V) = \dot \bigcup_{\gamma \in \Gamma \left/ \spitz{\gamma_0}\right.} \gamma \, \pi_{\gamma_0}^{- 1} \rund{\overline \pi^{\, - 1}(V) \cap U} \, .
\]

So similar to the case of a regular point one gets a whole isomorphism $\rund{\left.\overline \pi \right|_U , \rund{|_{\pi_{\gamma_0}^{- 1} \rund{\overline \pi^{\, - 1}(V) \cap U} }}_{V \subset \overline \pi(U) \text{ open } }}$ of ringed spaces from 
$\rund{U, \left. \widetilde \S \right|_U}$ to $\rund{\overline \pi (U) , \S|_{\overline \pi (U)} }$ , where

\begin{eqnarray*}
&& |_{\pi_{\gamma_0}^{- 1} \rund{\overline \pi^{\, - 1}(V) \cap U}} : \\
&& \phantom{12} \S(V) = \schweif{f \in \O\rund{\pi^{- 1}(V)} \otimes \P^\cz \phantom{12} \Upsilon\text{-invariant } } \\
&& \phantom{12} \rightarrow \widetilde \S\rund{\overline \pi^{\, - 1}(V) \cap U} \\
&& = \schweif{\left.f \in \O\rund{\pi_{\gamma_0}^{- 1} \rund{\overline \pi^{\, - 1}(V) \cap U}} \otimes \P^\cz \, \right| \, f\rund{\widetilde{\gamma_0} z} = f}
\end{eqnarray*}

is simply the restriction map, which is an isomorphism of unital $\P^\cz$-algebras.

\begin{defin}
$\Phi(0) = \Pi\rund{z_0} \in_{\P} \Upsilon \backslash H$ is called an elliptic point of $\Upsilon \backslash H$~. Its body is $\overline{z_0}$ , which is an elliptic point of $\Gamma \backslash H$ .
\end{defin}

For giving a local trivialization of $\F$ at $\overline{z_0}$ first of all identify the ringed spaces $\rund{B, \O_B \otimes \P^\cz}$ , $\rund{\spitz{\eta} \backslash B \, , \, \widetilde \O_B \otimes \P^\cz}$ , 
$\rund{\left. \spitz{g^{- 1} \gamma_0 g} \right\backslash H \, , \, \widetilde \O_H \otimes \P^\cz}$ and $\X$ via $\Phi$ locally at \\
$0 \mapsto 0 \mapsto i \mapsto \overline{z_0}$ and define the $\widetilde \O_H$- sheaf module $\E_H$ on $\left.\spitz{g^{- 1} \gamma_0 g} \right\backslash H$ by

\[
\E_H(V) := \schweif{\left.f \in \O\rund{\pi_{g^{- 1} \gamma_0 g}^{- 1}(V)} \, \right| \, f|_{g^{- 1} \gamma_0 g} = f}
\]

for all $V \subset \left.\spitz{g^{- 1} \gamma_0 g} \right\backslash H$ open. Now we show that locally at $\overline{z_0}$ we have $\S$- sheaf module isomorphisms

\[
\F \rightarrow \E_H \otimes \P^\cz \rightarrow \E_B \otimes \P^\cz \rightarrow \O_B \otimes \P^\cz \, ,
\]

where $\E_B$ is a suitable $\widetilde O_B$- sheaf module on $\spitz{\eta} \backslash B$ . Similar to $\Phi$ the first isomorphism is given by the restriction maps $|_{\pi_{\gamma_0}^{- 1} \rund{\overline \pi^{\, - 1}(V) \cap U}}$ , 
$V \subset \overline \pi(U)$ open, followed by $|_g$ . The second one is given by

\[
|_c : \E_H(V) \rightarrow \E_B\rund{c^{- 1}(V)} \, , \, f \mapsto f|_c := f\rund{c(w)} j(c, w)^k
\]

for all $V \subset \left.\spitz{g^{- 1} \gamma_0 g} \right\backslash H$ , where $w$ denotes the standard holomorphic coordinate on $B$ and $j(c, w) := \frac{1 - i}{2} \, \frac{1}{z - 1} \in \O(B)$ is chosen such that $j(c, w)^2 = c'$ . \\

If $2 | k$ then $\E_B$ is given by

\[
\E_B(V) := \schweif{f \in \O\rund{\pi_\eta^{- 1}(V)} \, \left| \, f(\eta w) \, \eta^{k / 2} = f\right.}
\]

and the last isomorphism by

\[
\E_B(V) \rightarrow \O\rund{\rund{\sqrt[n]{\phantom{1}} \,}^{- 1} (V)} \, , \, f \mapsto f\rund{\sqrt[n]{w} \,} {\sqrt[n]{w}}^{\, k / 2} \, w^{- \ceil{\frac{k}{2 n}} }
\]

for all $V \subset \spitz{\eta} \backslash B$ open. \\

If $2 \not| k$ and $- 1 \notin \Gamma$ then automatically $2 \not| n$ , and so there exists a unique $\eps \in \spitz{\eta}$ such that $\eps^2 = \eta$ and so $\ord \ \eps = n$ . Now $\E_B$ is given by

\[
\E_B(V) := \schweif{\left.f \in \O\rund{\pi_\eta^{- 1}(V)} \, \right| \, f(\eta w) \, \eps^k = f}
\]

and the last isomorphism by

\[
\E_B(V) \rightarrow \O\rund{\rund{\sqrt[n]{\phantom{1}} \,}^{- 1} (V)} \, , \, f \mapsto f\rund{\sqrt[n]{w} \,} {\sqrt[n]{w}}^{\, (k + n) / 2} \, w^{- \ceil{\frac{k + n}{2 n}} }
\]

for all $V \subset \spitz{\eta} \backslash B$ open.

\end{quote}

{\it At cusps of $\Gamma \backslash H$ }:

\begin{quote}
Let $\overline{z_0} \in \Gamma \left\backslash \partial_{\pz^1} H \right.$ , $z_0 \in \partial_{\pz^1} H$ , be a cusp of $\Gamma \backslash H$ , and let $g \in G$ such that $g \infty = z_0$ and 
$\overline{g^{- 1} \Gamma g} \cap \overline{P^\infty} = \spitz{\overline{g_0}}$ . Let \\
$\overline \pi : \left. \spitz{g g_0 g^{- 1}} \right\backslash H \cup \{z_0\} \rightarrow X$ denote the canonical projection. Then $\pi = \overline \pi \circ \pi_{g g_0 g^{- 1}}$ , and similar to the case of an elliptic point

\begin{eqnarray*}
&& \psi: B \mathop{\longrightarrow}^{\frac{\log}{2 \pi i}} \left. \spitz{g_0} \right\backslash H \cup \{\infty\} \mathop{\longrightarrow}^\id \left. \spitz{g_0} \right\backslash H \cup \{\infty\} \\
&& \phantom{1234} \mathop{\longrightarrow}^g \left. \spitz{g g_0 g^{- 1}} \right\backslash H \cup \{z_0\} \mathop{\longrightarrow}^{\overline \pi} X
\end{eqnarray*}

gives a locally biholomorphic map at $0 \mapsto \infty \mapsto \infty \mapsto z_0 \mapsto \overline{z_0}$ , so $\psi^{- 1}$ is a local chart of $X$ at $\overline{z_0}$ . \\

Now for giving a local $\P$-chart of $\X$ at $\overline{z_0}$ define the sheaf $\widetilde \O$ of complex algebras on $\left. \spitz{\overline{g_0}} \right\backslash H \cup \{\infty\}$ , the sheaf $\widetilde \S_\infty$ of $\P^\cz$-algebras on 
$\left. \spitz{\overline{g_0}} \right\backslash H \cup \{\infty\}$ and the sheaf $\widetilde \S_{z_0}$ of $\P^\cz$-algebras on $\left. \spitz{\overline{g g_0 g^{- 1}} } \right\backslash H \cup \schweif{z_0}$ by

\begin{eqnarray*}
&& \widetilde \O(V) := \left\{f \in \O\rund{\pi_{g_0}^{- 1}\rund{V \setminus \{\infty\}} } \, \right| \, f\rund{g_0 z} = f \text{ and } \\
&& \left. \phantom{1234567 \pi_{g_0}^{- 1}} f \text{ bounded at } \infty \text{ if } \infty \in V\right\}
\end{eqnarray*}

and

\begin{eqnarray*}
&& \widetilde \S_\infty(V) := \left\{f \in \O\rund{\pi_{g_0}^{- 1}\rund{V \setminus \{\infty\}} } \otimes \P^\cz \, \right| \, f\rund{\widetilde{g_0} z} = f \text{ and } \\
&& \left. \phantom{123456789 \pi_{g_0}^{- 1}} f \text{ bounded at } \infty \text{ if } \infty \in V\right\}
\end{eqnarray*}

for all $V \subset \left.\spitz{\overline{g_0}} \right\backslash H \cup \{\infty\}$ open, and Finally,

\begin{eqnarray*}
&& \widetilde \S_{z_0}(V) := \left\{f \in \O\rund{\pi_{g g_0 g^{- 1}}^{- 1}\rund{V \setminus \schweif{z_0}} } \otimes \P^\cz \, \right| \, f\rund{\widetilde{g_0} z} = f \\
&& \left. \phantom{1234 \pi_{g_0 g^{- 1}}^{- 1} \P^\cz} \text{ and } f \text{ bounded at } z_0 \text{ if } z_0 \in V\right\}
\end{eqnarray*}

for all $V \subset \left.\spitz{\overline{g g_0 g^{- 1}} } \right\backslash H \cup \schweif{z_0}$ open, where $\widetilde{g_0} \in g^{- 1} \Upsilon g$ is the unique element such that $\widetilde{g_0}^\# = g_0$ . Again we will see that $\psi$ extends to a local 
isomorphism

\begin{eqnarray*}
&& \Psi: \rund{B, \O_B \otimes \P^\cz} \rightarrow \rund{\left. \spitz{g_0} \right\backslash H \cup \{\infty\} \, , \, \widetilde \O \otimes \P^\cz} \\
&& \phantom{12} \rightarrow \rund{\left. \spitz{g_0} \right\backslash H \cup \{\infty\} \, , \, \widetilde \S_\infty} \\
&& \phantom{12} \rightarrow \rund{\left. \spitz{g g_0 g^{- 1}} \right\backslash H \cup \schweif{z_0} \, , \, \widetilde \S_{z_0}} \rightarrow \X
\end{eqnarray*}

of ringed spaces at $0 \mapsto \overline{z_0}$ , and therefore $\Psi^{- 1}$ is a local $\P$-chart of $\X$ at $\overline{z_0}$ . The first isomorphism is induced by $\frac{\log}{2 \pi i}$ . Now let $\chi$ and $\Omega$ again be given by lemma \ref{paramparab} 
taken $\widetilde{g_0}$ instead of $g$ . This leads to the commuting diagram

\[
\begin{array}{ccccc}
\phantom{12,} H & \mathop{\longrightarrow}\limits^{\Omega} & H \phantom{12,} \\
g_0 \downarrow & \circlearrowleft & \downarrow \widetilde{g_0} \\
\phantom{12,} H & \mathop{\longrightarrow}\limits_{\Omega} & H \phantom{12,}
\end{array}
\]

of $\P$-automorphisms inducing the second isomorphism, and the third isomorphism is induced by $g$ as an automorphism of $H$~. The last isomorphism is obtained by the same procedure as the one for elliptic points using an open neighbourhood \\
$U \subset \spitz{g g_0 g^{- 1}} \backslash H \cup \schweif{z_0}$ of $z_0$ such that

\[
\left.\overline \pi \right|_U : U \rightarrow \overline \pi (U) \mathop{\subset}\limits_{ \text{ open } } X
\]

is biholomorphic.

\begin{defin}
$\Psi(0) \in_{\P} \X$ is called a cusp of $\Upsilon \backslash H$ . Its body is $\overline{z_0} \in X$ , which is a cusp of $\Gamma \backslash H$ .
\end{defin}

Observe that in general there is no interpretation of $\Psi(0)$ as a $\P$-point of $\partial_{\pz^1} H$ . \\

For giving a local trivialization of $\F$ at $\overline{z_0}$ identify the ringed spaces $\rund{B, \O_B \otimes \P^\cz}$ , $\rund{\spitz{g_0} \backslash H \cup \{\infty\} \, , \, \widetilde \O \otimes \P^\cz}$ and $\X$ via $\Psi$ locally at 
$0 \mapsto \infty \mapsto \overline{z_0}$ . Now we show that locally at $\overline{z_0}$ we have $\S$- sheaf module isomorphisms

\[
\F \rightarrow \E \otimes \P^\cz \rightarrow \O_B \otimes \P^\cz \, ,
\]

where $\E$ is a suitable $\widetilde \O$- sheaf module on $\spitz{g_0} \backslash H$ . Similar to $\Psi$ the first isomorphism is given by the restriction maps $|_{\pi_{\gamma_0}^{- 1} \rund{\overline \pi^{\, - 1}(V) \cap U}}$ , 
$V \subset \overline \pi(U)$ open, followed by $|_g$ and $|_\Omega$ . \\

\begin{defin}
Assume $- 1 \notin \Gamma$ . Then either $g_0 \in g^{- 1} \Gamma g$ or $- g_0 \in g^{- 1} \Gamma g$ . So the cusp $\overline{z_0} \in X$ of $\Gamma \backslash H$ is called even (odd) iff $g_0 \in g^{- 1} \Gamma g$ 
(resp. $- g_0 \in g^{- 1} \Gamma g$~).
\end{defin}

If either $2 | k$ or $2 \not| k$ , $- 1 \notin \Gamma$ and $\overline{z_0}$ even then $\E$ is given by

\begin{eqnarray*}
&& \E(V) := \left\{f \in \O\rund{\pi_{g_0}^{- 1}(V \setminus \{\infty\}) } \, \right| \, f|_{g_0} = f \text{ and } \\
&& \left. \phantom{1234567 \pi_{g_0}^{- 1}} f \text{ bounded at } \infty \text{ if } \infty \in V\right\}
\end{eqnarray*}

and since $f|_{g_0}(z) = f(z + 1)$ for all $z \in H$ the last isomorphism by

\[
\E(V) \rightarrow \O\rund{\rund{\frac{\log}{2 \pi i}}^{- 1}(V)} \, , \, f \mapsto f\rund{\frac{\log w}{2 \pi i}}
\]

for all $V \subset \left.\spitz{g_0} \right\backslash H$ open, where $w$ denotes the standard holomorphic coordinate on $B$ . \\

If $2 \not| k$ , $- 1 \notin \Gamma$ and $\overline{z_0}$ odd then $\E$ is given by

\begin{eqnarray*}
&& \E(V) := \left\{f \in \O\rund{\pi_{g_0}^{- 1}(V \setminus \{\infty\}) } \, \right| \, f|_{- g_0} = f \text{ and } \\
&& \left. \phantom{1234567 \pi_{g_0}^{- 1}} f \text{ bounded at } \infty \text{ if } \infty \in V\right\}
\end{eqnarray*}

and since $f|_{- g_0}(z) = - f(z + 1)$ for all $z \in H$ the last isomorphism by

\[
\E(V) \rightarrow \O\rund{\rund{\frac{\log}{2 \pi i}}^{- 1}(V)} \, , \, f \mapsto f\rund{\frac{\log w}{2 \pi i}} e^{- \frac{1}{2} \log w}
\]

for all $V \subset \left.\spitz{g_0} \right\backslash H$ open.

\end{quote}

\section{$\mathcal{P}$-points of the {\sc Teichmüller} space} \label{Teich}

Let $g \in \nz$ and $\T_g$ be the {\sc Teichmüller} space for genus $g$ . Let us recall some basic properties. $\T_g$ is a complex domain of dimension

\[
N_g = \left\{\begin{array}{ll}
0 & \text{ if } g = 0 \\
1 & \text{ if } g = 1 \\
3 (g - 1) & \text{ if } g \geq 2
\end{array}\right.  \, .
\]

For every $\ba \in \T_g$ let $S(\ba)$ be its corresponding compact {\sc Riemann} surface of genus $g$ . Then all these compact {\sc Riemann} surfaces $S(\ba)$ , $\ba \in \T_g$ , glue together to a holomorphic family $\pi: \Xi_g \rightarrow \T_g$ over $\T_g$ 
with $S(\ba) := \pi^{- 1}(\ba)$~, $\ba \in \T_g$ , in particular $\pi$ is a holomorphic submersion, and the moduli space of compact {\sc Riemann} surfaces of genus $g$ is given by

\[
\M_g = \T_g \left/ \Gamma_g\right.
\]

with a certain discrete subgroup $\Gamma_g \sqsubset \Aut \T_g$ .

\begin{examples}
\end{examples}

\begin{itemize}
\item[(i)] $\T_0$ consists of one single point with the {\sc Riemann} sphere $\pz^1$ as corresponding compact {\sc Riemann} surface.
\item[(ii)] $\T_1 = H$ , $\Gamma_g = SL(2, \zz)$ , and $S(a) = \cz / (\zz + a \zz)$ for all $a \in H$ .
\end{itemize}

Now let $\ba \in \T_g$ be fixed and $S(\ba)$ be given by $U_1, \dots, U_n \subset \cz$ open together with the glueing data

\[
\sigma_{i j} : U_i \mathop{\supset}_{ \text{ open } } U_{i j} \rightarrow U_{j i} \mathop{\subset}_{ \text{ open } } U_j
\]

biholomorphic. Then of course $\rund{U_i}_{i = 1, \dots, n}$ is an open cover of $S(\ba)$ , and after some refinement of this open cover we may assume that \\
$H^1\rund{\rund{U_i}_{i \in \{1, \dots, n\}} , T X} \simeq H^1(X , T X) \simeq T_{\ba} \T_g$ canonically, and then using charts of $\Xi_g$ where the projection $\pi$ is just given by projecting onto the first $N_g$ coordinates we see that there exist an open 
neighbourhood $B$ of $\ba$ in $\T_g$ and families $\rund{U_{i j}^{(\bw)}}_{\bw \in B}$ , $\rund{\sigma_{i j}^{(\bw)}}_{\bw \in B}$ , $i, j = 1, \dots, n$ , such that

\begin{itemize}
\item[\{i\}]

\[
U_{i j}^{(\bw)} \subset U_i
\]

is open and

\[
\sigma_{i j}^{(\bw)} : U_{i j}^{(\bw)} \rightarrow U_{j i}^{(\bw)}
\]

is biholomorphic for all $i, j = 1, \dots, n$ and $\bw \in B$ ,

\item[\{ii\}]

\[
\dot \bigcup_{\bw \in B} U_{i j}^{(\bw)} \subset U_i \times B
\]

is open and

\[
\dot \bigcup_{\bw \in B} U_{i j}^{(\bw)} \rightarrow \cz \, , \, (z, \bw) \mapsto \sigma_{i j}^{(\bw)}(z)
\]

is holomorphic for all $i, j = 1, \dots, n$ ,

\item[\{iii\}] $U_{i j}^{(\ba)} = U_{i j}$ and $\sigma_{i j}^{(\ba)} = \sigma_{i j}$ for all $i, j = 1, \dots, n$ , and Finally,

\item[\{iv\}] $S(\bw)$ is given by the glueing data

\[
\sigma_{i j}^{(\bw)} : U_i \mathop{\supset}_{ \text{ open } } U_{i j}^{(\bw)} \rightarrow U_{j i}^{(\bw)} \mathop{\subset}_{ \text{ open } } U_j
\]

for all $\bw \in B$ , and so in particular we have a $\cz$-linear map

\begin{eqnarray*}
\Omega: T_{\ba} \T_g &\rightarrow& Z^1\rund{\rund{U_i}_{i \in \{1, \dots, n\}} , T X} \, , \\
\bv &\mapsto& \rund{\rund{\sigma_{i j}'}^{- 1} \rund{\left.\partial_t \sigma_{i j}^{(\ba + t \bv)} \right|_{t = 0}} }_{i, j \in\{1, \dots, n\}}
\end{eqnarray*}

such that $\eckig{\phantom{1}} \circ \Omega : T_{\ba} \T_g \rightarrow H^1(X , T X)$ is an isomorphism, where

\[
\eckig{\phantom{1}} : Z^1\rund{\rund{U_i}_{i \in \{1, \dots, n\}} , T X} \rightarrow H^1\rund{\rund{U_i}_{i \in \{1, \dots, n\}} , T X} \hookrightarrow H^1(X , T X)
\]

denotes the canonical projection.
\end{itemize}

Now let $\widetilde \ba \in_{\P} B$ such that ${\widetilde \ba}^\# = \ba$ . Then according to (vii) of section \ref{param} we can assign to $\widetilde \ba$ the $\P$- {\sc Riemann} surface $S\rund{\widetilde \ba} = \pi^{- 1}\rund{\widetilde \ba}$ . 
$S\rund{\widetilde a}$ is given by the local $\P$-charts $U_1, \dots, U_n \subset \cz$ with the $\P$- glueing data

\[
\sigma_{i j}^{\widetilde a} : U_i \mathop{\supset}_{ \text{ open } } U_{i j} \rightarrow_{\P} U_{j i} \mathop{\subset}_{ \text{ open } } U_j \, ,
\]

which are $\P$-isomorphisms, and clearly its body is $S(a)$ . The purpose of this section is to show that any $\P$- {\sc Riemann} surface $\X$ with compact body $X := \X^\#$ of genus $g$ can be realized as a $\P$-point of the {\sc Teichmüller} 
space $\T_g$ , which is of general interest. Before doing so we need a lemma:

\begin{lemma} \label{autoparaRiemann}
\item[(i)] Let $\X$ be a $\P$- {\sc Riemann} surface such that $X := \X^\#$ is compact of genus $g \geq 2$ , and let $\Phi$ be a $\P$-automorphism of $\X$ with $\Phi^\# = \id$ . Then $\Phi = \id$ .

\item[(ii)] Let $a \in_{\P} H = \T_1$ . Then the $\P$-automorphisms of $S(a) = \cz / (\zz + a \zz)$ with $\id$ as body are given by the translations 

\[
t_b : \cz / (\zz + a \zz) \rightarrow_{\P} \cz / (\zz + a \zz) \, , \, z \mapsto z + b \, ,
\]

$b \in \I^\cz$ , where $\I \lhd \P$ denotes the unique maximal ideal in $\P$ .
\end{lemma}

Recall that the $\P$- {\sc Riemann} surface $S(a) = \cz / (\zz + a \zz)$ in (ii) can be written as $\rund{\cz \left/ \rund{\zz + a^\# \zz}\right., \S}$ , where $\S$ is given by

\[
\S(U) := \schweif{\left.f \in \O\rund{\pi^{- 1}(U)} \otimes \P^\cz \, \right| \, f(z + m + n a) = f(z) \text{ for all } m, n \in \zz}
\]

for all $U \subset \cz \left/ \rund{\zz + a^\# \zz}\right.$ open and $\pi: \cz \rightarrow \cz \left/ \rund{\zz + a^\# \zz}\right.$ denotes the canonical projection. \\

{\it Proof:} via induction on $N$ , where $\I^N = 0$ , $\I$ being the unique maximal ideal of $\P$ . If $N = 1$ then $\P = \rz$ , and so both assertions are trivial. Now assume $\I^{N + 1} = 0$ , and define

\[
\Q := \P \left/ \I^N\right.
\]

with unique maximal ideal $\J := \I \left/ \I^N\right. \lhd \Q$ having $\J^N = 0$ , and let ${}^\natural : \P \rightarrow \Q$ be the canonical projection. \\

For proving (i) let $\X$ be a $\P$- {\sc Riemann} surface with body $X := \X^\#$ , compact of genus $\geq 2$ , and $\Phi$ be a $\P$-automorphism of $\X$ with $\Phi^\# = \id$ . Let $\X$ be given by $U_1, \dots, U_n \subset \cz$ open together with the 
$\P$- glueing data

\[
\varphi_{i j} : U_i \mathop{\supset}_{ \text{ open } } U_{i j} \rightarrow_{\P} U_{j i} \mathop{\subset}_{ \text{ open } } U_j \, .
\]

Then in the local $\P$-charts $U_i$ , $i = 1, \dots, n$ , $\Phi$ is given by $\P$-automorphisms $\Phi_i : U_i \rightarrow_{\P} U_i$ having $\Phi_i^\# = \id$ and

\[
\begin{array}{ccc}
\phantom{123456,} U_{i j} & \mathop{\longrightarrow}\limits^{\varphi_{i j}} & U_{j i} \phantom{12345,} \\
\left.\Phi_i \right|_{U_{i j}} \uparrow & \circlearrowleft & \uparrow \left.\Phi_j \right|_{U_{j i}} \\
\phantom{123456,} U_{i j} & \mathop{\longrightarrow}\limits_{\varphi_{i j}} & U_{j i} \phantom{12345,}
\end{array}  \, .
\]

By induction hypothesis $\Phi^\natural = \id$ , so $\Phi_i = \id + f_i$ with suitable $f_i \in \O\rund{U_{i j}} \otimes \rund{\I^\cz}^N$ , and a straight forward calculation using $\I^{N + 1} = 0$ shows that $f_i$ glue together to an element 
$f \in H^0(X, T X) \otimes \rund{\I^\cz}^N$ , but $H^0(X, T X) = 0$ since $X$ is of genus $\geq 2$ . So all $\Phi_i = \id$ , and so $\Phi = \id$ . \\

For proving (ii) let $a \in_{\P} H$ and $\Phi$ be a $\P$-automorphism of $S(a) = \cz / (\zz + a \zz)$~. By induction hypothesis we already know that there exists $c \in \I^\cz$ such that $\Phi^\natural = t_{c^\natural}$ . Then the lift of $\Phi$ as a 
$\P$-automorphism of $\cz$ commuting with translations by $\zz + a \zz$ is given by $z \mapsto z + c + f$ with a suitable $f \in \O(\cz) \otimes \rund{\I^\cz}^N$ . But then $f$ has to be invariant under translations by $\zz + a^\# \zz$ , and so it has to be 
constant $\in \rund{\I^\cz}^N$ . Define $b := c + f$ . $\Box$ \\

\begin{theorem} \label{paramTeich}
Let $\X$ be a $\P$- {\sc Riemann} surface, whose body $X = S(\ba)$ is a compact {\sc Riemann} surface of genus $g$ , $\ba \in \T_g$ . Then there exists a unique $\widetilde \ba \in_{\P} \T_g$ such that

\begin{itemize}
\item[(i)] $\widetilde \ba ^\# = \ba$ and
\item[(ii)] there exists a $\P$-isomorphism

\[
\Phi: S\rund{\widetilde \ba} \rightarrow_{\P} \X
\]

having $\Phi^\# = \id$ .
\end{itemize}
\end{theorem}

By lemma \ref{autoparaRiemann} $\Phi$ is uniquely determined by $\X$ if and only if $g \geq 2$ . \\

{\it Proof:} again via induction on $N$ , where $\I^N = 0$ , $\I$ being the unique maximal ideal of $\P$ . If $N = 1$ then $\P = \rz$ , and so the assertion is again trivial. So assume $\I^{N + 1} = 0$ , and let $\X$ be a $\P$- {\sc Riemann} surface with body 
$X := \X^\# = S(\ba)$ , $\ba \in \T_g$ . Let $\X$ be given by $U_1, \dots, U_n \subset \cz$ open together with the $\P$- glueing data

\[
\varphi_{i j} : U_i \mathop{\supset}_{ \text{ open } } U_{i j} \rightarrow_{\P} U_{j i} \mathop{\subset}_{ \text{ open } } U_j \, .
\]

Then $\rund{U_i}_{i \in \{1, \dots, n\}}$ forms an open cover of $X$ , and after maybe some refinement of this open cover we may again assume that

\[
H^1\rund{\rund{U_i}_{i = 1, \dots, n} , T X} \simeq H^1(X, T X) \simeq T_{\ba} \T_g
\]

canonically. So let $B$ be an open neighbourhood of $\ba$ in $\T_g$ and $\rund{U_{i j}^{(\bw)}}_{\bw \in B}$ and $\rund{\sigma_{i j}^{(\bw)}}_{\bw \in B}$ be families such that the conditions \{i\} - \{iv\} are fulfilled. Again define $\Q := \P \left/ \I^N\right.$ with 
unique maximal ideal $\J := \I \left/ \I^N\right. \lhd \Q$ having $\J^N = 0$ , and let ${}^\natural : \P \rightarrow \Q$ be the canonical projection. Let $X^\natural$ be the $\Q$- {\sc Riemann} surface given by the $\Q$- glueing data

\[
\varphi_{i j}^\natural : U_i \mathop{\supset}_{ \text{ open } } U_{i j} \rightarrow_{\P} U_{j i} \mathop{\subset}_{ \text{ open } } U_j \, .
\]

Then evidently $\rund{\X^\natural}^\# = X$ , and therefore by induction hypothesis there exists a unique $\bb \in_{\P} B$ with $\bb^\# = \ba$ such that there exists a $\Q$-isomorphism $\Psi: S\rund{\bb^\natural} \rightarrow_{\Q} \X^\natural$ having 
$\Psi^\# = \id$ . In the local $\Q$-charts $U_i$ , $i = 1, \dots, n$~, $\Psi$ is given by $\Q$-automorphisms $\Psi_i : U_i \rightarrow_{\Q} U_i$ having $\Psi^\# = \id$ and

\[
\begin{array}{ccc}
\phantom{123456,} U_{i j} & \mathop{\longrightarrow}\limits^{\varphi_{i j}^\natural} & U_{j i} \phantom{12345,} \\
\left.\Psi_i \right|_{U_{i j}} \uparrow & \circlearrowleft & \uparrow \left.\Psi_j \right|_{U_{j i}} \\
\phantom{123456,} U_{i j} & \mathop{\longrightarrow}\limits_{\sigma^{\rund{\bb^\natural}}_{i j}} & U_{j i} \phantom{12345,}
\end{array}  \, .
\]

Let $\widetilde{\Psi_i} : U_i \rightarrow_{\P} U_i$ , $i = 1, \dots, n$ , be arbitrary such that $\widetilde{\Psi_i}^\natural = \Psi_i$ , and so $\widetilde{\Psi_i}$ , $i = 1, \dots, n$ , are automatically $\P$-automorphisms, and define

\[
\rho_{i j} := \rund{\left.\widetilde{\Psi_j} \right|_{U_{j i}} }^{- 1} \circ \varphi_{i j} \circ \left.\widetilde{\Psi_i} \right|_{U_{i j}} : U_{i j} \rightarrow_{\P} U_{j i}  \, .
\]

Then for all $i, j = 1, \dots, n$ since $\rho_{i j}^\natural = \sigma_{i j}^{\rund{\bb^\natural}} = \rund{\sigma_{i j}^{(\bb)} }^\natural$ and $\sigma_{i j}^{(\bb)}$ is a $\P$-automorphism we see that

\[
\rho_{i j}(z) = \sigma_{i j}^{(\bb)}\rund{z + \omega_{i j}(z)}
\]

with some suitable $\omega_{i j} \in \O\rund{U_{i j}} \otimes \rund{\I^\cz}^N$ , $i, j \in \{1, \dots, n\}$ . Now an easy calculation shows that the $\omega_{i j}$ , $i, j \in \{1, \dots, n\}$ , form an element 
$\omega \in Z^1\rund{\rund{U_i}_{i \in \{1, \dots, n\}} , T X} \otimes \rund{\I^\cz}^N$ , and therefore there exist \\
$\bV \in T_\ba \T_g \otimes \rund{\I^\cz}^N$ and $\mu \in C^0\rund{\rund{U_i}_{i \in \{1, \dots, n\}} , T X} \otimes \rund{\I^\cz}^N$ such that $\omega = \Omega(\bV) + \delta \mu$ , where $\delta$ denotes the coboundary operator on the cocomplex associated to 
the sheaf cohomology of $T X$ and the open cover $X = \bigcup_{i = 1, \dots, n} U_i$  . In the local charts $U_i$ , $i = 1, \dots, n$ , $\mu$ is given by some $\mu_i \in \O\rund{U_i} \otimes \rund{\I^\cz}^N$ , $i = 1, \dots, n$ . Now define 
$\widetilde \ba := \bb + \bV$ and

\[
\Phi_i : U_i \rightarrow_{\P} U_i \, , \, z \mapsto \widetilde{\Psi_i}(z) - \mu_i(z)  \, ,
\]

$i = 1, \dots, n$ . Then ${\widetilde \ba}^\# = \ba$ , $\Phi_i^\# = \id$ , and so $\Phi_i$ , $i = 1, \dots, n$ , are automatically $\P$-automorphisms. A straight forward calculation shows that

\[
\begin{array}{ccc}
\phantom{123456,} U_{i j} & \mathop{\longrightarrow}\limits^{\varphi_{i j}} & U_{j i} \phantom{12345,} \\
\left.\Phi_i \right|_{U_{i j}} \uparrow & \circlearrowleft & \uparrow \left.\Phi_j \right|_{U_{j i}} \\
\phantom{123456,} U_{i j} & \mathop{\longrightarrow}\limits_{\sigma^{\rund{\widetilde \ba}}_{i j}} & U_{j i} \phantom{12345,}
\end{array}  \, ,
\]

and so $\Phi_i$ , $i = 1, \dots, n$ , glue together to a $\P$-isomorphism $\Phi: S\rund{\widetilde \ba} \rightarrow_{\P} \X$ with $\Phi^\# = \id$ . This proves the existence of $\widetilde \ba$ . \\

For proving uniqueness we may assume without loss of generality that $g \geq 1$ . Let also $\bc \in_{\P} B$ such that $\bc^\# = \ba$ , and let $\Lambda: S\rund{\widetilde \ba} \rightarrow_{\P} S\rund{\bc}$ be a $\P$-isomorphism such that 
$\Lambda^\# = \id$ . We will show that $\bc = \widetilde \ba$ . By induction hypothesis we already know that ${\widetilde \ba}^\natural = \bc^\natural$ , so $\bW := \bc - \widetilde \ba \in T_\ba \T_g \otimes \rund{\I^\cz}^N$ , and $\Lambda^\natural$ is a 
$\Q$-automorphism of $S\rund{\widetilde \ba}^\natural = S\rund{{\widetilde \ba}^\natural}$ . \\

{\it First case}: $g \geq 2$ .

\begin{quote}
Then by lemma \ref{autoparaRiemann} (ii) we directly know that $\Lambda^\natural = \id$ , and so in the local charts $U_i$ , $i = 1, \dots, n$ , $\Lambda$ is given by $\P$-automorphisms $\Lambda_i$ of $U_i$ with

\[
\begin{array}{ccc}
\phantom{123456,} U_{i j} & \mathop{\longrightarrow}\limits^{\sigma^{\rund{\bc}}_{i j}} & U_{j i} \phantom{12345,} \\
\left.\Lambda_i \right|_{U_{i j}} \uparrow & \circlearrowleft & \uparrow \left.\Lambda_j \right|_{U_{j i}} \\
\phantom{123456,} U_{i j} & \mathop{\longrightarrow}\limits_{\sigma^{\rund{\widetilde \ba}}_{i j}} & U_{j i} \phantom{12345,}
\end{array}  \, ,
\]

and $\Lambda_i^\natural = \id$ . Therefore $\Lambda_i = \id + f_i$ with suitable \\
$f_i \in \O\rund{U_i} \otimes \rund{\I^\cz}^N$ .
\end{quote}

{\it Second case}: $g = 1$ .

\begin{quote}
By lemma \ref{autoparaRiemann} (i) we know that $\Lambda^\natural = t_{d^\natural}$ with some $d \in \I^\cz$ . So

\[
\Lambda \circ t_{- d} : S\rund{\widetilde \ba} \rightarrow_{\P} S(\bc)
\]

is a $\P$-isomorphism having $\rund{\Lambda \circ t_{- d}}^\natural = \id$ . So again in the local charts $U_i$ , $i = 1, \dots, n$ , $\Lambda \circ t_{- d}$ is given by $\id + f_i$ where $f_i \in \O\rund{U_i} \otimes \rund{\I^\cz}^N$ .
\end{quote}

In both cases now an easy calculation shows that $\Omega({\bW}) = \delta f$ , \\
$f := \rund{f_i}_{i = 1, \dots, n}$ considered as an element of $C^0\rund{\rund{U_i}_{i = 1, \dots, n} , T X} \otimes \rund{\I^\cz}^N$~. Therefore $\overline \Omega (\bW) = \b0$ , and so $\bW = \b0$ . $\Box$ \\

\begin{cor} \label{paramTeich0}
Let $\X$ be a $\P$- {\sc Riemann} surface with body $\X^\# = \pz^1$ . Then there exists a $\P$-isomorphism $\Phi: \pz^1 \rightarrow_{\P} \X$ having $\Phi^\# = \id$ .
\end{cor}

\section{The main result} \label{sect main}

Now we return to the backbone of the article. So let $\X = (X, \S)$ , \\
$X := \Gamma \backslash H \cup \{ \text{ cusps of } \Gamma \backslash H\}$ , be the $\P$- {\sc Riemann } surface constructed in section \ref{quot} , and let $g$ be the genus of $X$ . Then by theorem 
\ref{paramTeich} we may identify $\X$ with $S\rund{\widetilde \ba }$ and so $X$ with $S(\ba)$ for some $\widetilde \ba \in_{\P} \T_g$ and $\ba := {\widetilde \ba}^\#$ . \\

Let $\widetilde{e_1}, \dots, \widetilde{e_R} \in_{\P} \X = S\rund{\widetilde \ba}$ , $\rho = 1, \dots, R$ , be the elliptic points and $\widetilde{s_1}, \dots, \widetilde{s_S} \in_{\P} \X$ , $\sigma = 1, \dots, S$ , the cusps of $\Upsilon \backslash H$ . Then \\
$e_\rho := \widetilde{e_\rho}^\# \in X = S(\ba)$ , $\rho = 1, \dots, R$ , $s_\sigma := \widetilde{s_\sigma}^\# \in X$ , $\sigma = 1, \dots, S$~, are automatically the elliptic points resp. cusps of $\Gamma \backslash H$ . Let $n_\rho \in \nz$ denote 
the period of the elliptic point $e_\rho$ , $\rho = 1, \dots, R$ . Let $U_\rho \subset X$ and $V_\sigma \subset X$ be pairwise disjoint open connected coordinate neighbourhoods of $e_\rho$ , $\rho = 1, \dots, R$ , and $s_\sigma$ , $\sigma = 1, \dots, S$ , 
resp.. Then via common local charts on $S(\bw)$ , $\bw \in W$ , $W \subset \T_g$ a suitably small open neighbourhood of $\ba$ , we can identify the sets $W \times U_\rho$ , $\rho = 1, \dots, R$ , $W \times V_\sigma$ , $\sigma = 1, \dots, S$ , with pairwise 
disjoint open sets of $\Xi_g$ such that 

\[
\begin{array}{ccc}
\phantom{123,} W \times U_\rho & \hookrightarrow & \Xi_g \phantom{1,} \\
\pr_1 \downarrow & \circlearrowleft &  \downarrow \pi \\
\phantom{123,} W & \subset & \T_g \phantom{1,}
\end{array}  \, ,
\]

$\rho = 1, \dots, R$ , and

\[
\begin{array}{ccc}
\phantom{123,} W \times V_\sigma & \hookrightarrow & \Xi_g \phantom{1,} \\
\pr_1 \downarrow & \circlearrowleft &  \downarrow \pi \\
\phantom{123,} W & \subset & \T_g \phantom{1,}
\end{array}  \, ,
\]

$\sigma = 1, \dots, S$ . This gives us at the same time embeddings $U_\rho , V_\sigma \hookrightarrow S(\bw)$ , $\rho = 1, \dots, R$ , $\sigma = 1, \dots S$ , $\bw \in W$ , as pairwise disjoint open sets. Now let 

\[
\widetilde E := \rund{\widetilde{e_1}, \dots, \widetilde{e_R}, \widetilde{s_1}, \dots, \widetilde{s_S}} \in_{\P} \X^{R + S}
\]

and

\[
\U := U_1 \times \cdots \times U_R \times V_1 \times \cdots \times V_S \subset X^{R + S} \, ,
\]

which is an open neighbourhood of \\
$E := \widetilde E^\# = \rund{e_1, \dots, e_R, s_1, \dots, s_S} \in X^{R + S}$ , and at the same time for all $\bw \in W$ it is identified with some open set in $S(\bw)^{R + S}$ . Let

\[
\left.\Xi_g \right|_W := \pi^{- 1}(W) = \dot \bigcup_{\bw \in W} S(\bw)  \, .
\]

Then $(\pi, \id) : \left.\Xi_g \right|_W \times \U \rightarrow W \times \U$ is a familiy of compact {\sc Riemann} surfaces $S(\bw) \times \{\bu\} = S(\bw)$ , $(\bw, \bu) \in W \times \U$ . Furthermore, let

\[
U_0 := \dot \bigcup_{(\bw, \bu) \in W \times \U} \rund{S(\bw) \setminus \schweif{u_1, \dots, u_R, v_1, \dots, v_S}} \times \{\bu\} \subset \left.\Xi_g \right|_W \times \U
\]

open and dense, where $\bu = \rund{u_1, \dots, u_R, v_1, \dots, v_S}$ . Then

\[
\left.\Xi_g \right|_W \times \U = U_0 \cup \bigcup_{\rho = 1}^R \rund{W \times U_\rho \times \U} \cup \bigcup_{\sigma = 1}^S \rund{W \times V_\sigma \times \U}
\]

is a finite open cover, and we can define holomorphic line bundles on $\left.\Xi_g \right|_W \times \U$ via trivializations on $U_0$ and on each $W \times U_\rho \times \U$~, $\rho = 1, \dots, R$ , and $W \times V_\sigma \times \U$~, 
$\sigma = 1, \dots, S$ , and transition functions $\phi_{0 \rho} \in \O\rund{U_0 \cap \rund{W \times U_\rho \times \U}}$~, $\rho = 1, \dots, R$ , resp. $\psi_{0 \sigma} \in \O\rund{U_0 \cap \rund{W \times V_\rho \times \U}}$ , $\sigma = 1, \dots, S$ . \\

Furthermore, let $T^*_{\rel}$ denote the relative cotangent bundle of the family $\left.\Xi_g \right|_W \times \U \rightarrow W \times \U$ of compact {\sc Riemann} surfaces, see for example section 10.1 of \cite{Schlich} . It is a holomorphic line bundle on 
$\left.\Xi_g \right|_W \times \U$ such that $\left.T^*_{\rel} \right|_{S(\bw) \times \{\bu\}}$ is the cotangent bundle of $S(\bw)$ for all $(\bw, \bu) \in W \times \U$~. Therefore $\left.T^*_{\rel} \right|_{S(\bw) \times \{\bu\}} = T^* S(\bw)$ even for all $\P$-points 
$(\bw, \bu) \in_{\P} W \times \U$ . \\

Finally, let $\F$ denote the holomorphic $\P$- line bundle over $\X = S\rund{\widetilde \ba}$ defined in section \ref{quot} having $M_k(\Upsilon) = H^0(\F)$ . \\

{\it First we treat the case $2 | k$ .} \\

We define the holomorphic line bundle $L_k \rightarrow \left.\Xi_g \right|_W \times \U$ by the transition functions $\phi_{0 \rho}(z) := \rund{z - u_\rho}^{- k / 2 + \ceil{\frac{k}{2 n_\rho}} }$ , $\rho = 1, \dots, R$ , and 
$\psi_{0 \sigma}(z) := \rund{z - v_\sigma}^{- k / 2}$ , $\sigma = 1, \dots, S$ , where $z$ denotes a local coordinate on $U_\rho$ resp. $V_{\sigma}$ . Then for each $(\bw, \bu) \in W \times \U$ the holomorphic sections of $L_k|_{S(\bw) \times \{\bu\}}$ are 
precisely the meromorphic functions on $S(\bw)$ with poles at the points $u_\rho \in U_\rho \hookrightarrow S(\bw)$ of order at most $k / 2 - \ceil{\frac{k}{2 n_\rho}}$ , $\rho = 1, \dots, R$~, and poles at the points $v_\sigma \in V_\sigma \hookrightarrow 
S(\bw)$ , $\sigma = 1, \dots, S$ , of order at most $k / 2$ and holomorphic at all other points of $S(\bw)$ . \\

Furthermore, we define the holomorphic line bundle $C \rightarrow \left.\Xi_g \right|_W \times \U$ by the transition functions $\phi_{0 \rho}(z) := 1$ , $\rho = 1, \dots, R$ , and $\psi_{0 \sigma}(z) := z - v_\sigma$ , $\sigma = 1, \dots, S$ . Then for each 
$(\bw, \bu) \in W \times \U$ the holomorphic sections of $C|_{S(\bw) \times \{\bu\}}$ are precisely the holomorphic functions on $S(\bw)$ vanishing at $v_\sigma \in V_\sigma \hookrightarrow S(\bw)$ , $\sigma = 1, \dots, S$ . Therefore clearly 
$\deg C|_{S(\bw) \times \{\bu\}} = - S$ for all $(\bw, \bu) \in W \times \U$ , and we identify the sections of $C$ with ordinary holomorphic functions on $\left.\Xi_g \right|_W \times \U$ vanishing on $\schweif{z = v_\sigma}$ , $\sigma = 1, \dots, S$ . \\

Finally, we define the line bundle $L'_k := {T^*_{\rel}}^{\otimes (k / 2)} \otimes L_k$ over $\left.\Xi_g \right|_W \times \U$ .

\begin{lemma} \label{isoMkSk}

We have isomorphisms $\F \simeq \left.L'_k \right|_{S\rund{\widetilde \ba} \times \schweif{\widetilde E}}$ and so \\
$M_k(\Upsilon) \simeq H^0\rund{\left.L'_k \right|_{S\rund{\widetilde \ba} \times \schweif{\widetilde E}} }$ , the last isomorphism mapping $S_k(\Upsilon)$ precisely to 
$H^0\rund{\left.\rund{L'_k \otimes C}\right|_{S\rund{\widetilde \ba} \times \schweif{\widetilde E}} }$ .

\end{lemma}

{\it Proof:} Since $g' = j(g, z)^2$ for all $g \in_{\P} G$ regarded as a $\P$-automorphism of $H$ , identifying the trivial and the cotangent bundle on $H$ we see that 
$\F \simeq \rund{T^* \X}^{\otimes (k / 2)} = \left.{T_{\rel}^*}^{\otimes (k / 2)}\right|_{S\rund{\widetilde \ba} \times \schweif{\widetilde E}}$ on $X' := (\Gamma \backslash H) \setminus \schweif{e_1, \dots, e_r}$ . \\

Now let $\Phi_1 , \dots, \Phi_R$ and $\Psi_1 , \dots, \Psi_S$ denote the local $\P$-charts of $\X$ at $e_1, \dots, e_R$ and $s_1, \dots, s_S$ resp. given in section \ref{quot} . Then via these local $\P$-charts the elliptic points 
$\widetilde e_1, \dots, \widetilde e_R$ and the cusps $\widetilde s_1, \dots, \widetilde s_S \in_{\P} \X$ of $\Upsilon \backslash H$ are identified with the ordinary point $0 \in B$ . \\

A straight forward calculation shows that the holomorphic sections of $\rund{T^* \X}^{\otimes (k / 2)}$ regarded as holomorphic sections of $\F$ on $X'$ vanish in the local $\P$-chart $\Phi_\rho$ at $0$ of order $k / 2 - \ceil{\frac{k}{2 n}}$ for all 
$\rho = 1, \dots, R$ and in the local $\P$-chart $\Psi_\sigma$ at $0$ of order $k / 2$ for all $\sigma = 1, \dots, S$ . This proves the first statement. \\

Furthermore, if $f \in M_k(\Upsilon) = H^0(\F)$ one sees that $f \in S_k(\Upsilon)$ iff $f$ vanishes in the local $\P$-chart $\Psi_\sigma$ at $0$ for all $\sigma = 1, \dots, S$ , which proves the second statement. $\Box$ \\

Observe that 

\begin{eqnarray} \label{degree}
\deg \left. L'_k \right|_{S(\bw) \times \{\bu\}} &=& \frac{k}{2} \deg \left. T^*_{\rel} \right|_{S(\bw) \times \{\bu\}} + \deg \left. L_k \right|_{S(\bw) \times \{\bu\}}  \notag \\
&=& k (g - 1) + \frac{k}{2} (R + S) - \sum_{\rho = 1}^R \ceil{\frac{k}{2 n_\rho}}
\end{eqnarray}

is independent of the point $(\bw, \bu) \in W \times \U$ .

\begin{lemma} \label{globsectMkSk}

\[
\dot \bigcup_{(\bw, \bu) \in W \times \U} H^0\rund{\left. L'_k \right|_{S(\bw) \times \{\bu\}} }
\]

is a holomorphic vector bundle over $W \times \U$ , containing

\[
\dot \bigcup_{(\bw, \bu) \in W \times \U} H^0\rund{\left.\rund{L'_k \otimes C}\right|_{S(\bw) \times \{\bu\}} }
\]

as a holomorphic sub vector bundle.
\end{lemma}

{\it Proof:} We will use theorem 5 in section 10.5 of \cite{GrauRem} , which says the following:

\begin{quote}
If $\dim_\cz H^i\rund{X_y, \underline V_y}$ is independent of $y \in Y$ then all sheaves $f_{(i)}\rund{\underline V}$ are locally free and all maps 

\[
f_{y, i}: f_{(i)}\rund{\underline V} \left/ \m_y f_{(i)}\rund{\underline V}\right. \rightarrow H^i\rund{X_y, \underline V_y}
\]

are isomorphisms.
\end{quote}

Hereby $f: X \rightarrow Y$ denotes a holomorphic family of compact complex manifolds $X_y := f^{- 1}(y)$ , $y \in Y$ , $\underline V$ a holomorphic vector bundle over $X$ and $f_{(i)}\rund{\underline V}$ , $i \in \nz$ , the higher direct image sheaves of 
$\underline V$ under $f$ . $\underline V_y := \left. \underline V \right|_{X_y}$~, $\m_y \lhd \O_Y$ denotes the maximal ideal of holomorphic functions on $Y$ vanishing at the point $y \in Y$ , and finally, $f_{y, i}: f_{(i)}\rund{\underline V} \left/ \m_y f_{(i)}\rund{\underline V}\right. \rightarrow H^i\rund{X_y, \underline V_y}$ denotes the canonical homomorphism. Recall that the sheaf $f_{(0)}\rund{\underline V}$ is given by the assignment 
$U \mapsto H^0\rund{f^{- 1}(U), \underline V}$ for all $U \subset Y$ open.

So we have to show that $\dim H^0\rund{\left. L'_k \right|_{S(\bw) \times \{\bu\}} }$ and $\dim H^0\rund{\left.\rund{L'_k \otimes C}\right|_{S(\bw) \times \{\bu\}} }$ are independent of the point $(\bw, \bu) \in W \times \U$~. We use formula (\ref{degree}) . The case 
$g = 0$ is trivial.

\begin{itemize}

\item[Let $g = 1$] . Then $T^*_{\rel}$ is the trivial bundle, and $R + S \geq 1$ .

\begin{itemize}
\item[for $k = 0$ :] $L_k'$ is trivial. If $S = 0$ then $L_k' \otimes C$ is trivial, if $S \geq 1$ then $\deg \left. \rund{L_k' \otimes C} \right|_{S(\bw) \times \{\bu\}} < 0$ so $H^0\rund{\left. \rund{L_k' \otimes C} \right|_{S(\bw) \times \{\bu\}} } = \{0\}$ .
\item[for $k = 2$ :] If $S = 0$ then $L_k'$ is trivial, if $S \geq 1$ then $\deg \left. L_k' \right|_{S(\bw) \times \{\bu\}} = S \geq 1$ and so $\dim H^0\rund{\left. L'_k \right|_{S(\bw) \times \{\bu\}} } = \deg \left. L_k' \right|_{S(\bw) \times \{\bu\}} = S$ . 

$L_k' \otimes C$ is trivial.
\item[for $k \geq 4$ :] $\deg \left. L_k' \right|_{S(\bw) \times \{\bu\}} , \deg \left. \rund{L_k' \otimes C} \right|_{S(\bw) \times \{\bu\}} \geq R + S \geq 1$ .
\end{itemize}

\item[Let $g \geq 2$] .

\begin{itemize}
\item[for $k = 0$ :] $L_k'$ is trivial. If $S = 0$ then $L_k' \otimes C$ is trivial, if $S \geq 1$ then $\deg \left. \rund{L_k' \otimes C} \right|_{S(\bw) \times \{\bu\}} < 0$ .
\item[for $k = 2$ :] If $S = 0$ then $L_k' = T^*_{\rel}$ so $\dim H^0\rund{\left. L'_k \right|_{S(\bw) \times \{\bu\}} } = g$ , if $S \geq 1$ then $\deg \left. L_k' \right|_{S(\bw) \times \{\bu\}} = 2 (g - 1) + S \geq 2 g - 1$ and so 
$\dim H^0\rund{\left. L'_k \right|_{S(\bw) \times \{\bu\}} } = \deg \left. L_k' \right|_{S(\bw) \times \{\bu\}} - g + 1 = g - 1 + S$ .

$L_k' \otimes C = T^*_{\rel}$ .
\item[for $k \geq 4$ :] $\deg \left. L_k' \right|_{S(\bw) \times \{\bu\}} , \deg \left. \rund{L_k' \otimes C} \right|_{S(\bw) \times \{\bu\}} \geq 2 g - 1$ . $\Box$
\end{itemize}

\end{itemize}

\begin{theorem}[main theorem] \label{main}

We have isomorphisms

\[
\begin{array}{ccc}
S_k(\Upsilon) & \simeq & S_k(\Gamma) \otimes \P^\cz \\
\cap & \circlearrowleft & \cap \\
M_k(\Upsilon) & \simeq & M_k(\Gamma) \otimes \P^\cz \\
{}_{{}^\#} \searrow & \circlearrowleft & \swarrow_{\id \otimes {}^\#} \\
 & M_k(\Gamma) &
\end{array} \, .
\]

\end{theorem}

{\it Proof:} By lemma \ref{globsectMkSk} after maybe replacing $W$ and $U_\rho$ , $\rho = 1, \dots, R$~, $V_\sigma$ , $\sigma~=~1, \dots, S$~, by smaller open neighbourhoods of $\ba \in \T_g$ resp. $e_\rho , s_\sigma \in X$ 
there exist $F_1, \dots, F_r \in H^0\rund{L'_k \otimes C} \hookrightarrow H^0\rund{L'_k}$ and $F_{r + 1} , \dots, F_{r'} \in~H^0\rund{L'_k}$ such that $\rund{\left.F_\rho \right|_{S(\bw) \times \{\bu\}} }_{\rho \in \{1, \dots, r\}}$ is a basis of 
the $\cz$-vectorspace $H^0\rund{\left.\rund{L'_k \otimes C}\right|_{S(\bw) \times \{\bu\}} }$ and $\rund{\left.F_\rho \right|_{S(\bw) \times \{\bu\}} }_{\rho \in \{1, \dots, r'\}}$ is a basis of the $\cz$-vectorspace $H^0\rund{\left. L'_k \right|_{S(\bw) \times \{\bu\}} }$ 
for all $(\bw, \bu) \in W \times \U$ . Now define 

\[
f_\rho := \left.F_\rho \right|_{S\rund{\widetilde \ba} \times \schweif{\widetilde E}} \in H^0\rund{\left. L'_k \right|_{S\rund{\widetilde \ba} \times \schweif{\widetilde E}} } \, ,
\]

$\rho = 1, \dots, r'$ . Then

\[
f_1, \dots, f_r \in H^0\rund{\left.\rund{L'_k \otimes C}\right|_{S\rund{\widetilde \ba} \times \schweif{\widetilde E}} } \, ,
\]

$\rund{f_1^\#, \dots, f_r^\#}$ is a basis of $H^0\rund{\rund{\left.\rund{L'_k \otimes C}\right|_{S\rund{\widetilde \ba} \times \schweif{\widetilde E}} }^\#}$ , and $\rund{f_1^\#, \dots, f_{r'}^\#}$ is a basis of 
$H^0\rund{\rund{\left. L'_k \right|_{S\rund{\widetilde \ba} \times \schweif{\widetilde E}} }^\#}$ .

\[
\rund{\left.\rund{L'_k \otimes C}\right|_{S\rund{\widetilde \ba} \times \schweif{\widetilde E}} }^\# = \left.\rund{L'_k \otimes C}\right|_{S(\ba) \times \{E\}}
\]

and

\[
\rund{\left. L'_k \right|_{S\rund{\widetilde \ba} \times \schweif{\widetilde E}} }^\# = \left. L'_k \right|_{S(\ba) \times \{E\}} \, .
\]

One obtains the result combining this with lemmas \ref{globsect} and \ref{isoMkSk} . $\Box$ \\

{\it Now we treat the case $2 \not| k$ and $- 1 \notin \Gamma$ .} \\

Let $s_1, \dots, s_{S'} \in X$ be the even and $s_{S' + 1}, \dots, s_S \in X$ the odd cusps of $\Gamma \backslash H$ . \\

We define the holomorphic line bundle $L_k \rightarrow \left.\Xi_g \right|_W \times \U$ by the transition functions $\phi_{0 \rho}(z) := \rund{z - u_\rho}^{- k - 1 + 2 \ceil{\frac{k+ n_\rho}{2 n_\rho}} }$ , $\rho = 1, \dots, R$ , 
$\psi_{0 \sigma}(z) := \rund{z - v_\sigma}^{- k}$ , $\sigma = 1, \dots, S'$ and $\psi_{0 \sigma}(z) := \rund{z - v_\sigma}^{- k + 1}$ , $\sigma = S' + 1, \dots, S$ . \\

Furthermore, we define the holomorphic line bundle $C \rightarrow \left.\Xi_g \right|_W \times \U$ by the transition functions $\phi_{0 \rho}(z) := 1$ , $\rho = 1, \dots, R$ , $\psi_{0 \sigma}(z) := z - v_\sigma$~, $\sigma = 1, \dots, S'$ , and 
$\psi_{0 \sigma}(z) := 1$ , $\sigma = S' + 1, \dots, S$ . Clearly $\deg C|_{S(\bw) \times \{\bu\}} = - S'$ for all $(\bw, \bu) \in W \times \U$ , and we identify the sections of $C$ with ordinary holomorphic functions on $\left.\Xi_g \right|_W \times \U$ vanishing on 
$\schweif{z = v_\sigma}$ , $\sigma = 1, \dots, S'$ . \\

Now let $F := \F^\#$ , which is an ordinary holomorphic line bundle on $X =~S(\ba)$ . Then obviously $M_k(\Gamma) = H^0(F)$ , and a straight forward calculation similar to the proof of lemma \ref{isoMkSk} shows that \\
$F^{\otimes 2} \simeq \left.\rund{{T^*_\rel}^{\otimes k} \otimes L_k} \right|_{S(\ba) \times \{E\}}$ . So after maybe replacing $W$ and $U_\rho$ , $\rho = 1, \dots, R$ , $V_\sigma$ , $\sigma = 1, \dots, S$ , by smaller open neighbourhoods of $\ba \in \T_g$ 
resp. $e_\rho , s_\sigma \in X$ we may assume that there exists a unique line bundle $L'_k \rightarrow W \times \U$ such that $F \simeq \left.L'_k \right|_{S(\ba) \times \{E\}}$ and ${L'_k}^{\otimes 2} = {T^*_\rel}^{\otimes k} \otimes L_k$ .

\begin{lemma} \label{isoMkSk2}

We have isomorphisms $\F \simeq \left.L'_k \right|_{S\rund{\widetilde \ba} \times \schweif{\widetilde E}}$ and so \\
$M_k(\Upsilon) \simeq H^0\rund{\left.L'_k \right|_{S\rund{\widetilde \ba} \times \schweif{\widetilde E}} }$ , the last isomorphism mapping $S_k(\Upsilon)$ to $H^0\rund{\left.\rund{L'_k \otimes C}\right|_{S\rund{\widetilde \ba} \times \schweif{\widetilde E}} }$ .

\end{lemma}

{\it Proof:} By the same method as in the proof of lemma \ref{isoMkSk} one shows that

\[
\F^{\otimes 2} \simeq \left.\rund{{T^*_\rel}^{\otimes k} \otimes L_k} \right|_{S(\widetilde \ba) \times \schweif{\widetilde E}} = \rund{\left.L'_k \right|_{S\rund{\widetilde \ba} \times \schweif{\widetilde E}} }^{\otimes 2} \, . 
\]

Since $\F^\# = F \simeq \left.L'_k \right|_{S(\ba) \times \{E\}} = \rund{\left.L'_k \right|_{S\rund{\widetilde \ba} \times \schweif{\widetilde E}} }^\#$ the first assertion is a trivial consequence of the lemma \ref{root} below. \\

For proving the last statement one just has to observe that if $f \in M_k(\Upsilon) = H^0(\F)$ then $f \in S_k(\Upsilon)$ iff $f$ vanishes in the local $\P$-chart $\Psi_\sigma$ at $0$ for all $\sigma = 1, \dots, S'$ , where $\Psi$ denotes the local 
$\P$-chart of $\X$ at the cusp $s_\sigma$ , $\sigma = 1, \dots, S'$ . $\Box$ \\

\begin{lemma} \label{root}
Let $\E$ be a holomorphic $\P$- line bundle over the complex $\P$-manifold $\M = (M, \S)$ , $n \in \nz$ and $F$ be a holomorphic line bundle over $M$ such that $F^{\otimes n} = \E^\#$ . Then there exists an up to isomorphism unique $\P$- line bundle $\F$ 
over $\M$ such that $\F^\# = F$ and $\F^{\otimes n} \simeq \E$ .
\end{lemma}

{\it Proof:} Let $\E$ and $F$ be given by the local trivializations on $U_i \subset M$ open, $i \in I$ , $M = \bigcup_{i \in I} U_i$ , with $\P$- transition functions $\varphi_{i j} \in \S\rund{U_i \cap U_j}$ resp. transition functions 
$\psi_{i j} \in \O\rund{U_i \cap U_j}$ , $i, j \in I$ . Without loss of generality we may even assume that $\S|_{U_i} \simeq \O_{U_i} \otimes \P$ for all $i \in I$ . \\

Since $\cz \setminus \{0\} \rightarrow \cz \setminus \{0\} \, , \, z \mapsto z^n$ is locally biholomorphic there exist unique $\widetilde{\psi_{i j}} \in \S\rund{U_i \cap U_j} = \O\rund{U_i \cap U_j} \otimes \P^\cz$ such that $\widetilde{\psi_{i j}}^\# = \psi_{i j}$ and 
$\psi_{i j}^n = \varphi_{i j}$ , $i, j \in I$ . Now for proving existence define $\F$ via the local trivializations $\S|_{U_i}$ together with transition functions $\widetilde{\psi_{i j}} \in \S\rund{U_i \cap U_j}$ . For proving uniqueness let $\F'$ be another 
holomorphic $\P$- line bundle on $\M$ with body $F$ and ${\F'}^{\otimes n} \simeq \E$ . After maybe some refinement of the open cover $M = \bigcup_{i \in I} U_i$ we may assume without loss of generality that also $\F'$ admits local trivializations 
$\S|_{U_i}$ together with $\P$-transition functions $\eps_{i j} \in \S\rund{U_i \cap U_j}$ such that $\eps_{i j}^\# = \psi_{i j}$ and $\eps_{i j}^n = \varphi_{i j}$ , $i, j \in I$ . Therefore $\eps_{i j} = \widetilde{\psi_{i j}}$ , $i, j \in I$ , and so $\F' \simeq \F$ . 
$\Box$ \\

Again

\begin{eqnarray} \label{degree2}
\deg \left. L'_k \right|_{S(\bw) \times \{\bu\}} &=& \frac{k}{2} \deg \left. T^*_{\rel} \right|_{S(\bw) \times \{\bu\}} + \frac{1}{2} \deg \left. L_k \right|_{S(\bw) \times \{\bu\}} \notag \\
&=& k (g - 1) + \frac{k + 1}{2} R - \sum_{\rho = 1}^R \ceil{\frac{k + n_\rho}{2 n_\rho}} \notag \\
&& + \frac{k - 1}{2} S + \frac{S'}{2}
\end{eqnarray}

is independent of the point $(\bw, \bu) \in W \times \U$ , so automatically $2 \, | \, S'$ .

\begin{lemma} \label{globsectMkSk2}
Let $g \leq 1$ or $k \geq 3$ . Then 

\[
\dot \bigcup_{(\bw, \bu) \in W \times \U} H^0\rund{\left. L'_k \right|_{S(\bw) \times \{\bu\}} }
\]

is a holomorphic vector bundle over $W \times \U$ containing

\[
\dot \bigcup_{(\bw, \bu) \in W \times \U} H^0\rund{\left.\rund{L'_k \otimes C}\right|_{S(\bw) \times \{\bu\}} }
\]

as a holomorphic sub vector bundle.
\end{lemma}

{\it Proof:} same as the proof of lemma \ref{globsectMkSk} in the case $2 | k$ , now using formula (\ref{degree2}) .

\begin{itemize}

\item[Let $g = 1$] .

\begin{itemize}
\item[for $k = 1$ :] If $S' = 0$ then $L_k' = L_k' \otimes C$ and $L_k'^{\otimes 2}$ is trivial. So either $L_k'$ is trivial or $H^0\rund{\left. L_k' \right|_{S(\bw) \times \{\bu\}} } = \{0\}$ . If $S' \geq 2$ then \\
$\deg \left. L_k' \right|_{S(\bw) \times \{\bu\}} = \frac{S'}{2} \geq 1$ and \\
$\deg \left. \rund{L_k' \otimes C} \right|_{S(\bw) \times \{\bu\}} = - \frac{S'}{2} < 0$ .
\item[for $k \geq 3$ :] If $k = 3$ , $S = 0$ and all $n_\rho = 2$ then  \\
$L_k' = L_k' \otimes C$ and $L_k'^{\otimes 2}$ is trivial. In all other cases $\deg \left. L_k' \right|_{S(\bw) \times \{\bu\}} , \deg \left. \rund{L_k' \otimes C} \right|_{S(\bw) \times \{\bu\}} \geq 1$ .
\end{itemize}

\item[Let $g \geq 2$] .

\begin{itemize}
\item[for $k \geq 3$ :] $\deg \left. L_k' \right|_{S(\bw) \times \{\bu\}} , \deg \left. \rund{L_k' \otimes C} \right|_{S(\bw) \times \{\bu\}} \geq 3 (g - 1) \geq 2 g - 1$  .  $\Box$
\end{itemize}

\end{itemize}

\begin{theorem}[main theorem] \label{main2}

If $g \leq 1$ or $k \geq 3$ then we have isomorphisms

\[
\begin{array}{ccc}
S_k(\Upsilon) & \simeq & S_k(\Gamma) \otimes \P^\cz \\
\cap & \circlearrowleft & \cap \\
M_k(\Upsilon) & \simeq & M_k(\Gamma) \otimes \P^\cz \\
{}_{{}^\#} \searrow & \circlearrowleft & \swarrow_{\id \otimes {}^\#} \\
 & M_k(\Gamma) &
\end{array} \, .
\]

\end{theorem}

{\it Proof:} similar to the case $2 | k$ now using lemmas \ref{globsect} , \ref{isoMkSk2} and \ref{globsectMkSk2} . $\Box$ \\

In the case $g \geq 2$ and $k = 1$ the theorem indeed fails to be true. Here a counterexample: \\

Let $U \subset \cz$ be an open connected neighbourhood of $0$ , $\Xi \rightarrow U$ be a holomorphic family of compact {\sc Riemann} surfaces $S(w)$ of genus $g \geq 1$ and $L \rightarrow \Xi$ be a holomorphic 
line bundle (so automatically $\deg L|_{S(w)}$ is independent of the point $w$ ) having $\dim H^0\rund{L|_{S(w)} } < \dim H^0\rund{L|_{S(0)} }$ for all $w \in U \setminus \{0\}$ (so automatically $0 \leq \deg L|_{S(w)} \leq 2 (g - 1)$ ). Then there exists 
$f \in H^0\rund{L|_{S(0)} }$ with the following property:

\begin{quote}
$f$ admits {\bf no} extension to 'compact Riemann surfaces nearby'~, which means there exists no pair $\rund{U', F}$ , where $U' \subset U$ is an open neighbourhood of $0$ and $F \in H^0\rund{L|_{U'}}$ such that $F|_{S(w)} =~f$~.
\end{quote}

Now let $z_0 \in S(0)$ be arbitrary. After maybe replacing $U$ by a smaller open neighbourhood of $0$ we may fix a common local chart of all $S(w)$ , $w \in U$ , being a local coordinate neighbourhood of $S(0)$ at $z_0$ . Via this common local chart we 
may regard $z_0$ as a point of $S(w)$ for each $w \in U$~. Let $R \rightarrow \Xi$ be the holomorphic line bundle such that for each $w \in U$ the holomorphic sections of $R|_{S(w)}$ are the meromorphic functions on $S(w)$ which are holomorphic on 
$S(w) \setminus \{z_0\}$ and have a pole at $z_0 \in S(w)$ of order at most $d := 2 g - 1 - \deg L|_{S(w)}$ . So the holomorphic functions on $S(w)$ can be regarded as holomorphic sections of $R$ vanishing at $z_0$ of order at least $d$ . \\

Clearly $\deg (L \otimes R)|_{S(w)} = 2 g - 1$ , and so $\dim H^0\rund{\rund{L \otimes R}|_{S(w)} } = g$ is independent of the point $w \in U$ . So again by theorem 5 in section 10.5 of \cite{GrauRem} we see that

\[
\dot \bigcup H^0\rund{\rund{L \otimes R}|_{S(w)} }
\]

is a vector bundle of rank $g$ over $U$ . After maybe replacing $U$ by a smaller open neighbourhood of $0$ we can assume that there exists a frame $\rund{F_1, \dots, F_g} \in H^0(L \otimes R)^{\oplus g}$ . So there exists $F \in H^0(L \otimes R)$ having 
$F|_{S(0)} = f$ , and of course

\[
\Phi(F)(w) := \rund{\begin{array}{c}
F\rund{w, z_0} \\
\vdots \\
\partial_z^{d - 1} F\rund{w, z_0}
\end{array}} \in \O(U)^{\oplus d}
\]

must have an isolated zero at $w = 0$ for all such extensions $F$ of $f$ .

\begin{lemma} \label{examplezero}
There exists $N \in \nz$ such that $\ord_0 \, \Phi(F) \leq N$ for all extensions $F \in H^0(L \otimes R)$ of $f$ . 
\end{lemma}

{\it Proof:} Let $F \in H^0(L \otimes R)$ having $F|_{S(0)} = f$ be given. Using the language of germs at $0 \in U$ all other extensions of $f$ in 
$H^0(L \otimes R)$ are given by $\widetilde F = F + w \sum_{j = 1}^g \varphi_j F_j$ with $\varphi_j \in \O_{U, 0}$ , $j = 1, \dots, g$ .

\begin{eqnarray*}
\Phi\rund{\widetilde F}(w) &=& \rund{\begin{array}{c}
F\rund{w, z_0} \\
\vdots \\
\partial_z^{d - 1} F\rund{w, z_0}
\end{array}} \\
&& + w \rund{\begin{array}{ccc}
F_1\rund{w, z_0} & \dots & F_g\rund{w, z_0} \\
\vdots & & \\
\partial_z^{d - 1} F_1\rund{w, z_0} & \dots & \partial_z^{d - 1} F_g\rund{w, z_0}
\end{array}} \rund{\begin{array}{c}
\varphi_1 \\
\vdots \\
\varphi_g
\end{array}} \, .
\end{eqnarray*}

After trigonalization of the matrix using total pivot search with respect to the order at $w = 0$ one obtains up to multiplication with some element of 
$GL\rund{d, \O_{U, 0}}$ and permutation of the $\varphi_j$'s

\[
\Phi\rund{\widetilde F}(w) = \rund{\begin{array}{c}
H_1 \\
\vdots \\
H_d
\end{array}} + w \rund{\begin{array}{c}
\begin{array}{cccc}
a_1 & & & \\
 & \ddots & & b_{i j} \\
0 & & a_{d'} & 
\end{array} \\ \hline
0
\end{array}} \rund{\begin{array}{c}
\varphi_1 \\
\vdots \\
\varphi_g
\end{array}} \, ,
\]

$d' \in \{0, \dots, \min(d, g)\}$ , $H_i, a_i, b_{i j} \in \O_{U, 0}$ , $\ord_0 \, b_{i j} \geq \ord_0 \, a_i$ , $i = 1, \dots, d'$~, $j = i + 1, \dots, g$ . \\

Since the linear system $\Phi\rund{\widetilde F} = 0$ has no solution $\rund{\varphi_1, \dots, \varphi_g} \in \O_{U, 0}^{\oplus g}$ there exists at least one $i \in \{1, \dots, g\}$ such that $H_i \not= 0$ if $i \geq d' + 1$ resp. $\ord_0 \, H_i \leq \ord_0 \, a_i$ if 
$i \leq d'$ . So take $N := \ord_0 \, H_i$ . $\Box$ \\

Now let $\P := \rz[X] \left/ \rund{X^{N + 1} = 0}\right.$ and $\widetilde w := \overline X \in \P$ . Then $\widetilde w \in_{\P} U$ with $\widetilde w^\# = 0$ .

\begin{lemma}
There exists {\bf no} $\widetilde f \in H^0\rund{L|_{S\rund{\widetilde w}} }$ such that $\widetilde f^\# = f$ .
\end{lemma}

{\it Proof:} Assume $\widetilde f \in H^0\rund{L|_{S\rund{\widetilde w}} }$ having $\widetilde f^\# = f$ . Then $\widetilde f$ can be regarded as an element of $H^0\rund{(L \otimes R)|_{S\rund{\widetilde w}} }$ having a zero of order at least $d$ at 
$z_0 \in_{\P} S\rund{\widetilde w}$ . So there exist $\overline{a_1}, \dots, \overline{a_g} \in \P^\cz$ , $a_1, \dots, a_g \in \cz[X]$ , such that

\[
\widetilde f = \sum_{j = 1}^g \overline{a_j} \left. F_j \right|_{S\rund{\widetilde w}} \, .
\]

Since $f = \widetilde f^\# = \sum_{j = 1}^g \overline{a_j}^\# \left. F_j \right|_{S(0)}$ , we see that \\
$F :=  \sum_{j = 1}^g a_j(w) F_j \in H^0(L \otimes R)$ is an extension of $f$ . Finally, since 
$F|_{\widetilde w} = \widetilde f \in H^0\rund{(L \otimes R)|_{S\rund{\widetilde w}} }$ vanishes at $z_0 \in_{\P} S\rund{\widetilde w}$ of order at least $d$ we see that $\Phi(F)\rund{\widetilde w} = 0$ , and so $\Phi(F)$ has a zero at $w = 0$ of order at least 
$N + 1$ , which is a contradiction to lemma \ref{examplezero} . $\Box$ \\

Apply this to a suitable family of theta characteristics $L \rightarrow \Xi$ , $g \geq 2$~, which means a holomorphic line bundle $L$ having $L^{\otimes 2} = T^*_{\rel}$ , see \cite{Arb} , appendix B~. $\dim H^0\rund{L|_{S(w)} }$ is only constant 
$\mod 2$ . Since any holomorphic family of compact {\sc Riemann} surfaces of genus $g$ can be written as the pullback of a holomorphic map into $\T_g$ , and $\T_g$ , $g \geq 2$~, can be written as the moduli space of certain cocompact lattices in $G$ , 
see \cite{Nat} and example \ref{paramlatt} (ii) , we know that there exists a smooth map $\varphi: U \rightarrow G^{2 g} \, , \, w \mapsto \rund{A_1(w), B_1(w), \dots, A_g(w), B_g(w)}$ such that

\begin{itemize}
\item[(i)] all $A_1(w), B_1(w), \dots, A_g(w), B_g(w) \in G$ are hyperbolic generating a cocompact lattice $\Gamma_w$ without elliptic elements and $- 1 \notin \Gamma_w$ ,
\item[(ii)] $S(w) = \Gamma_w \backslash H$ ,
\item[(iii)] $L|_{S(w)}$ is obtained by the identification $(\gamma z, S) \sim \rund{z, j(\gamma, z) S}$ , \\
$\gamma \in \Gamma_w$ , in $H \times \cz$ .
\end{itemize}

In particular (iii) is guaranteed by lemma 11.1 in \cite{Nat} . So $\Upsilon := \Gamma_{\widetilde w}$ is a $\P$-lattice of $G$ with body $\Gamma := \Gamma_0$ , and by lemma \ref{isoMkSk2} we have \\
$M_1(\Upsilon) = S_1(\Upsilon) \simeq H^0\rund{L|_{S\rund{\widetilde w}} }$ .

\section{Body $SL(2, \zz)$ } \label{SL}

Let $\Upsilon$ be a $\P$-lattice of $G$ with body $\Gamma = SL(2, \zz)$ . Then of course this special case is easier to handle than the general theory. First of all $M_k(\Upsilon) = 0$ if $2 \not| \, k$ . So we can restrict our investigations the case $2 | k$ . \\

Let $\X = (X, \S)$ , $X := SL(2, \zz) \backslash H \cup \{\overline \infty\}$ , be the $\P$- {\sc Riemann} surface given by the construction of section \ref{quot} . Since $X$ is of genus $g = 0$ we may identify $X \simeq \pz^1$ , and from corollary \ref{paramTeich0} we 
know that there exists a $\P$-isomorphism $\Phi : X \rightarrow_{\P} \X$ with $\Phi^\# = \id$ . Let $\widetilde{e_1} , \widetilde{e_2} \in_{\P} \X$ be the two elliptic points of $\Upsilon \backslash H$ with bodies $e_1 := \overline i$ resp. 
$e_2 := \overline{e^{\frac{2}{3} \pi i}} \in X$ , which are precisely the elliptic points of $SL(2, \zz) \backslash H$ , and let $\widetilde s \in_{\P} \X$ be the cusp of $\Upsilon \backslash H$ with body $s := \overline \infty \in X$ , which is the cusp of 
$SL(2, \zz) \backslash H$ . One knows that $\Aut(X) \simeq SL(2, \cz) / \{\pm 1\}$ .

\begin{lemma}
There exists a unique $g \in_{\P} SL(2, \cz)$ such that $g^\# = 1$ , \\
$g e_1 = \Phi^{- 1}\rund{\widetilde{e_1}}$ , $g e_2 = \Phi^{- 1}\rund{\widetilde{e_2}}$ and $g s = \Phi^{- 1}\rund{\widetilde s}$ .
\end{lemma}

{\it Proof:} A straight forward computation shows that

\[
G \rightarrow X^3 \, , \, g \mapsto \rund{g e_1, g e_2, g s}
\]

is locally biholomorphic at $1$ , which proves the lemma. $\Box$ \\

Now a simple calculation using the local $\P$-charts of $\X$ given in section \ref{quot} shows that $\Phi \circ g$ uniquely lifts to a $\P$-automorphism $\Omega: H \rightarrow H$ having $\Omega^\# = \id$ such that 

\[
\begin{array}{ccc}
\phantom{1,} H & \mathop{\longrightarrow}\limits^{\Omega} & H \phantom{12} \\
\pi \downarrow & \circlearrowleft & \downarrow \Pi \\
\phantom{1,} X & \mathop{\longrightarrow}\limits_{\Phi \circ g} & \X \phantom{12}
\end{array} \, ,
\]

so automatically for all $\gamma \in \Upsilon$

\[
\begin{array}{ccc}
\phantom{123} H & \mathop{\longrightarrow}\limits^{\Omega} & H \phantom{1,} \\
\gamma^\# \downarrow & \circlearrowleft & \downarrow \gamma \\
\phantom{123} H & \mathop{\longrightarrow}\limits_{\Omega} & H \phantom{1,}
\end{array} \, .
\]

\begin{theorem}

For all $k \in 2 \nz$

\[
\Psi_k : M_k(\Upsilon) \rightarrow M_k(SL(2, \zz)) \otimes \P^\cz \, , \, f \mapsto f|_\Omega
\]

is a $\P^\cz$-module isomorphism mapping $S_k(\Upsilon)$ to $S_k(SL(2, \zz)) \otimes \P^\cz$ , and all $\Psi_k$ , $k \in 2 \nz$ , glue together to an isomorphism of $2 \, \nz$-graded $\P^\cz$-algebras

\[
\begin{array}{ccc}
\bigoplus\limits_{k \in 2 \nz} M_k(\Upsilon) & \rightarrow & \bigoplus\limits_{k \in 2 \nz} M_k(SL(2, \zz)) \otimes \P^\cz \\
\phantom{12345} {}_{{}^\#} \searrow & \circlearrowleft & \swarrow_{\id \otimes {}^\#} \phantom{1234567} \\
 & \bigoplus\limits_{k \in 2 \nz} M_k(SL(2, \zz)) & 
\end{array} \, .
\]

\end{theorem}

{\it Proof:} We use the notation of section \ref{quot} and the even case of section \ref{sect main} . Therefore we have to identify $X = X \times \schweif{\widetilde E}$ with $\X$ via $\Phi$ . $\Phi \circ g$ induces an identification 
$\rund{T^* \X}^{\otimes (k / 2)} \simeq \rund{T^* X}^{\otimes (k / 2)}$ , which restricted to $SL(2, \zz) \backslash H \subset X$ is given by $f \mapsto f|_{\Omega}$ , and $g$ as $\P$-automorphism of $X$ induces identifications 
$\left.L_k \right|_{X \times \schweif{\widetilde E}} \mathop{\rightarrow}\limits^{\sim} \left.L_k \right|_{X \times \{E\}}$ and $C|_{X \times \schweif{\widetilde E}} \mathop{\rightarrow}\limits^{\sim} C|_{X \times \{E\}}$ as holomorphic $\P$- line 
bundles given by $f \mapsto f(g z)$ , where we use a local coordinate $z$ on $X \simeq \pz^1$~. Since $\F \simeq \left.\rund{{T^*_{\rel}}^{\otimes (k / 2)} \otimes L_k} \right|_{X \times \schweif{\widetilde E}}$ by lemma \ref{isoMkSk} with body
$\F \simeq \left.\rund{{T^*_{\rel}}^{\otimes (k / 2)} \otimes L_k} \right|_{X \times \{E\}}$ and $\left. T^*_{\rel} \right|_{X \times \schweif{\widetilde E}} = \left. T^*_{\rel} \right|_{X \times \{E\}} = T^* X$ the claim follows. $\Box$ \\

\begin{small}

{\it Acknowledgement:} I have to thank M. {\sc Schlichenmaier} for many helpful comments during the writing procedure and the fonds national de recherche Luxembourg for funding my research stay at Luxembourg university.

\end{small}


\begin{thebibliography}{99}

\bibitem{Arb} {\sc Arbarello}, E. e. a., {\it Geometry of Algebraic Curves Vol. I}, Springer, Berlin New York, 1985.
\bibitem{GarlRagh} {\sc Garland}, H. and {\sc Raghunathan}, M. S., {\it Fundamental domains in ($\rz$-)rank 1 semisimple Lie groups}, Ann. Math. {\bf 92} (2) (1970), 279 - 326.
\bibitem{GrauRem} {\sc Grauert}, H. and {\sc Remmert}, R., {\it Coherent Analytic Sheaves. A Series of Comprehensive Studies in Mathematics}, Springer, Berlin Heidelberg, 1984.
\bibitem{Helga} {\sc Helgason}, S., {\it Differential Geometry and Symmetric Spaces}, Academic Press, New York, 1962.
\bibitem{Nat} {\sc Natanzon}, S. M., {\it Moduli of Riemann surfaces, Hurwitz-type spaces, and their super analogues}, Russ. Math. Surv. {\bf 54}:1 (1999), 61 - 117.
\bibitem{Ragh} {\sc Raghunathan}, M. S., {\it Discrete Subgroups of Lie Groups}, Springer, Berlin Heidelberg, 1972.
\bibitem{Schlich} {\sc Schlichenmaier}, M., {\it An Introduction to Riemann Surfaces, Algebraic Curves and Moduli Spaces, Theoretical and Mathematical Physics}, Springer, Berlin Heidelberg, 2007.


\end{thebibliography}
\end{document}